\documentclass[a4paper,10pt]{article}

\usepackage[a4paper,left=3cm,right=3cm,top=2.5cm,bottom=2.5cm]{geometry}
\usepackage{hyperref}
\usepackage{subfig}
\usepackage{pgfplots}
\usepackage{pgfplotstable}
\usepackage{mathtools}
\usepackage{multicol}
\usepackage{comment}
\usepackage{booktabs}
\pgfplotsset{compat=1.5}
\usepackage{amssymb}
\usepackage{url}
\usepackage{bm}

\usepackage{pdfsync}
\usepackage{float}
\usepackage{tabularx}
\usepackage{enumerate}
\usepackage{array}
\usepackage{xspace}
\usepackage{siunitx}
\usepackage{authblk}

\usepackage[draft,inline,marginclue]{fixme}
\usepackage{mathrsfs}
\usepackage{color, colortbl}
\usepackage{multirow}

\usepackage{amsthm}
\usepackage{amsmath}
\usepackage{stmaryrd}
\usepackage[ruled,vlined,linesnumbered]{algorithm2e}
\usepackage{graphicx}
\usepackage{booktabs}
\usepackage{enumitem}
\usepackage{cleveref}

\newcommand{\urom}{\ensuremath{\mathbf{u}_\text{POD}}}
\newcommand{\ufom}{\ensuremath{\mathbf{u}}}

\DeclareMathOperator*{\argmin}{\arg\!\min}
\newcolumntype{C}[1]{>{\centering\arraybackslash}m{#1}}

\definecolor{Gray}{gray}{0.9}

\newcommand{\RA}[1]{{\color{black}#1}}
\newcommand{\RB}[1]{{\color{black}#1}}

\begin{document}

\title{A DeepONet multi-fidelity approach for residual learning in reduced
  order modeling}

\author[1]{Nicola~Demo\footnote{nicola.demo@sissa.it}}
\author[1,2]{Marco~Tezzele\footnote{marco.tezzele@sissa.it}}
\author[1]{Gianluigi~Rozza\footnote{gianluigi.rozza@sissa.it}}

\affil[1]{Mathematics Area, mathLab, SISSA, via Bonomea 265, I-34136 Trieste, Italy}
\affil[2]{Oden Institute for Computational Engineering and
    Sciences, University of Texas at Austin, Austin, 78712, TX, United States}

\maketitle

\begin{abstract}
In the present work, we introduce a novel approach to enhance the
	precision of reduced order models by exploiting a multi-fidelity
	perspective and DeepONets. Reduced models provide a real-time numerical
	approximation by simplifying the original model. The error
	introduced by the such operation is usually neglected and sacrificed 
	in order to reach a fast computation. We propose to couple the
	model reduction to a machine learning residual learning, such
	that the above-mentioned error can be learned by a neural network
	and inferred for new predictions. We emphasize that the framework
	maximizes the exploitation of high-fidelity information,
	using it for building the reduced order model and for
	learning the residual. In this work, we explore the integration
	of proper orthogonal decomposition (POD), and gappy POD for sensors data,
	with the recent DeepONet architecture. Numerical investigations
	for a parametric benchmark function and a nonlinear parametric Navier-Stokes
	problem are presented.
\end{abstract}

\tableofcontents

\section{Introduction}
\label{sec:intro}

Multi-fidelity (MF) methods emerged as a solution to deal with complex
models, which usually need a high computational budget to be
solved~\cite{peherstorfer2018survey}. Such a framework aims to exploit not
only the so-called high-fidelity information, but also the response of
low-fidelity models in order to increase the accuracy of the
prediction.
This feature plays a fundamental role, especially for outer loop
applications such as uncertainty propagation and optimization, since
it allows to achieve good accuracy without 
requiring evaluating the high-fidelity model (typically expensive) at
every iteration. Thus, its employment is widespread for optimization
purposes, and among all the contributions in literature, we highlight the
successful application to naval engineering
problems~\cite{bonfiglio2018improving,bonfiglio2018multi,tezzele2021multi}, to multiple
fidelities modeling~\cite{forrester2007multi}, and in the presence of 
uncertainty~\cite{ng2014multifidelity}. 
All these cases, as well as many
others, build the correlation between the different fidelities by
involving Gaussian process regression (GPR).
Another approach
with nonlinear autoregressive schemes is described
in~\cite{perdikaris2017nonlinear,raissi2017inferring}, whereas
in~\cite{romor2021multi} a possible extension for high-dimensional
parameter spaces is investigated.
Recently, an alternative to such a probabilistic framework is offered by
neural networks, where the mapping between the low-fidelity model and the
high-fidelity one is learned by the network during the training procedure~\cite{zhang2021multi,meng2020composite,guo2021multi}. Among the
different types of architecture, DeepONet~\cite{lu2021learning,lin2023b} has been
proposed to approximate operators and it has been successfully applied to
MF problems in~\cite{lu2022multifidelity,howard2022multifidelity}. It
has also been successfully used to create a fast PDE-constrained
optimization method in~\cite{wang2021fast}. Another type of architecture that has been successfully applied to multi-fidelity data is the Bayesian neural network~\cite{meng2021multi}, resulting in a framework robust to noisy measurements.
We also highlight the employment of multi-fidelity techniques for uncertainty quantification. We cite~\cite{hart2022hyper, hart2022cal}
for a Bayesian framework capable to deal with model discrepancy using
different fidelities, whereas we refer to~\cite{farcas2022context} for an analysis of the trade-off between high- and low-fidelity data in a Monte Carlo estimation.

Reduced order modeling (ROM)~\cite{beohparour17, chinestaenc2017,RozzaStabileBallarin2022}
is a family of methods that aims at reducing the computational burden of evaluating complex 
models. Instead of combining data from heterogeneous models, ROM builds a
simplified model, typically from some high-fidelity information. Also, in this case, the
capabilities of ROM led to its diffusion in several industrial
contexts~\cite{morelli2021numerical, tezzele2021multi,
tezzele2020enhancing}, especially for optimization
tasks~\cite{benner2014model, amsallem2015design, zahr2015progressive, tezzele2018dimension, demo2020asga,
demo2021hull} or inverse problems~\cite{ghattas2021learning, ivagnes2022towards}.
In the ROM community, proper orthogonal decomposition (POD) is
one of the most employed methods to build the reduced
model~\cite{pichi2022aroma,qian2020lift,degruyter1,degruyter2,degruyter3}. Given
a limited set of high-fidelity data, POD is able to compute the reduced space of an
arbitrary rank which optimally (in a least squares sense) represents the 
data. In the last years, its diffusion led to several variants including
shifted POD~\cite{PapapiccoDemoGirfoglioStabileRozza2021,reiss2018shifted},
weighted POD~\cite{carere2021weighted,venturi2019weighted}, and gappy
POD~\cite{everson1995karhunen, willcox2006unsteady,
bui2004aerodynamic, mainini2015surrogate}. This latter exploit a compressive
sensing approach~\cite{bright2013compressive, brunton2014compressive,
kutz2017leveraging},  in order to use only a few information at
certain locations of the domain (sensors) to compute the approximation.
A generalization of gappy POD can be found in~\cite{adrian1975role},
where linear stochastic estimation allows the reconstruction of the
linear map between the available data and the system state by an $l_2$
minimization. 
A novel approach where such a relation between sensor data and the reduced state
is approximated in a nonlinear way employing neural networks can be
found in~\cite{nair2020leveraging}. 

In the present contribution, we explore the possibility of coupling these
two methodologies, MF and ROM, to enhance the accuracy of the
model. ROM indeed creates a simplified model from a few
high-fidelity data. Such approximation can be considered the low-fidelity model,
because of the projection error introduced by the ROM. In this context, MF could be adopted
in order to find the correlation between the original model and the ROM
one, resulting in a more precise prediction. We can therefore exploit
twice the collected high-fidelity data: initially, it is used to build
the reduced model, then again during the computation of the MF relation. From
this point of view, the proposed improvement does not need any additional
high-fidelity evaluations. Here we take into consideration the POD with interpolation or the gappy POD as
low-fidelity modeling techniques and the DeepONet to learn the
residual. POD with
interpolation~\cite{wang2012comparative, tezzele2022aroma09, gadalla2021les, demo2019cras} is
applied here for a completely data-driven approach, while gappy POD is
used in order to make the pipeline applicable even for sensor
data. The framework aims then to exploit the capability of POD models
for linear prediction, adding the nonlinear term through the DeepONet,
which can be viewed as a data-driven closure
model. See~\cite{xie2018data} for another data-driven modelling
approach to close ROMs, while for other recent works
that propose nonlinear model order reduction, we
cite~\cite{amsallem2012nonlinear, alla2017nonlinear,
kramer2019nonlinear, san2018neural, geelen2022operator, meneghetti2022aroma20, little2023nonlinear}.

The manuscript is organized as follows. In Section~\ref{sec:methods} we present the end-to-end numerical pipeline,
with a focus on POD, gappy POD, and DeepONets. We
continue in Section~\ref{sec:results} by showing the numerical
experiments, and finally we conclude with
Section~\ref{sec:conclusions} by summarizing the results and drawing
some future extensions.

\section{Methods}
\label{sec:methods}

This section is devoted to present the numerical methods used within the proposed approximation scheme, together with the methods used for comparison.
We describe their integration in order to provide a global overview, then we discuss in the following sections the algorithmic details.

Proper orthogonal decomposition (POD) is a widespread technique providing a linear model order reduction, particularly suited to deal with parametric problems~\cite{RozzaStabileBallarin2022,manzoni2016dimensionality,cueto2014model}.
Such a representation is computationally very cheap to acquire,
however, it suffers from the linear limitations of POD that may decrease its
accuracy, especially when dealing with nonlinear problems.

We are interested in efficiently computing a 
parametric field $\ufom (\mu)$ with $\ufom: P \to \mathcal V$, where $P$ is the parametric space, $\mathcal V$ a generic norm equipped vector space with $\text{dim}(\mathcal V) = n$.
POD-based ROMs compute the approximation $\urom(\mu)$
such that:
\begin{equation}
  \urom(\mu) \approx \ufom(\mu) = \urom(\mu) + \mathbf{r}(\mu),
  \label{eq:pod-mf}
\end{equation}
where $\mathbf r: P \to \mathcal V$ is the projection error introduced by the model order reduction, which we assume here to be dependent on the parameter. 
In a classical POD framework, this
residual $\mathbf{r}$ is usually neglected, due to its marginal contribution.
In the present contribution, we aim instead to
learn it by means of machine learning techniques, in order to improve the accuracy of the final prediction. Artificial neural networks (ANNs) can be used to model it, thanks to their general approximation capabilities, learning it by exploiting the snapshots already pre-computed to build the ROM.
In particular, dealing with parametric problems, we
exploit the DeepONet architecture to learn the
residual. The light computational demand to infer the DeepONet enables
a nonlinear but still real-time improvement of the POD model, at
the cost of additional training during the offline phase.

The only input needed by the proposed methodology is the numerical solutions
database $\{\mu_i, \ufom(\mu_i)\}_{i=1}^N$ computed by sampling the 
parameter space and exploiting any consolidated discretization method (e.g. finite element or finite volume method). These snapshots are combined in order to
find the POD space, which can be used for intrusive or non-intrusive ROM. We explore in this contribution only the non-intrusive (data-driven) approach, while future works will study the application to POD-Galerkin contexts. We investigate two options for the non-intrusive ROM:
\begin{itemize}
\item POD with radial basis functions (POD-RBF) interpolation, which enables the prediction of new solutions (for new parameters) by means of the above-mentioned interpolation technique. In this case, the ROM takes as input the actual parameter providing as output the approximated solution.
\item Gappy POD, which allows us to compute the approximated solution by providing only some sensor data thanks to a compressing strategy.
\end{itemize}
Once the ROM is built, we can exploit it to compute the low-fidelity representation of the original snapshots by passing the corresponding parameters (or sensor data).
The high-fidelity and low-fidelity databases are then used to learn the difference between them through the DeepONet network with the final aim of generalizing such residual even to unseen parameters and improving the final prediction. It is important to note that typically the space $\mathcal V$ is obtained by discretizing a generic $\mathbb R^d$ space. Depending on the complexity of the equation to solve and on the target accuracy, this kind of space can exhibit a high number of degrees of freedom. 

\RA{Approximating the error over such a high-dimensional space with a neural
network leads to two major issues: \textit{i}.) the number of the neurons in
the last layer is equal to the number of degrees of freedom of the space $\mathcal
V$, resulting in a model too large to treat; \textit{ii}.) the
parameter-to-error relation becomes too complex to be efficiently learned. Thus we extract the spatial coordinates of the degrees of freedom of
$\mathcal V$. Since we know the error (the difference between the original snapshots and the POD predictions) in any of these coordinates, we can
arrange the data in the format $\{(x_i, \mu_j, \mathbf r (\mu_j)_i)\,|\,
x_i \in \mathcal V \subset \mathbb R^d, \mu_j \in P, \mathbf r(\mu_j)_i \in \mathbb R\}$ where $i =
1, \dots, n$ and $j=1,\dots, N$, to isolate the spatial and parametric
dependency of the error. We can use such a dataset to learn the} \RB{scalar
error $r_{\text{net}}: \mathbb R^d \times P \to \mathbb R$} given the parametric and spatial coordinates.
\RA{In this way, the network maintains a limited number of output dimensions,
improving the identification of spatial recurrent patterns}. The loss function which is minimized
during the training procedure is then:
\begin{equation}
\mathcal L =  \frac{1}{N}\sum_{i=1}^N\| \mathbf r_\text{NN}(\mu_i) + \mathbf u_{\text{POD}}( \mu_i) - \mathbf u(\mu_i)\|^2,
\end{equation}
where $\mathbf r_\text{NN} (\mu)$ is not the high-dimensional output of a single network evaluation, but the array containing the result of the network inference for any spatial coordinates belonging to the discretized space such that \RB{$\mathbf r_\text{NN} (\mu) \equiv \begin{bmatrix} r_{\text{net}}(x_1, \mu) & r_{\text{net}}(x_2, \mu) & \dots & r_{\text{net}}(x_n, \mu) \end{bmatrix}$}, where $x_i \in \mathbb R^d$ for $i = 1, \dots, n$.
In the case of gappy POD, it is important to note that the DeepONet takes as input the sensors data and not the actual parameters.
We emphasize that the DeepONet training does not need any additional high-fidelity solutions besides those already collected for the POD space construction.

\begin{figure}[htb]
  \centering
  \includegraphics[width=1.\textwidth]{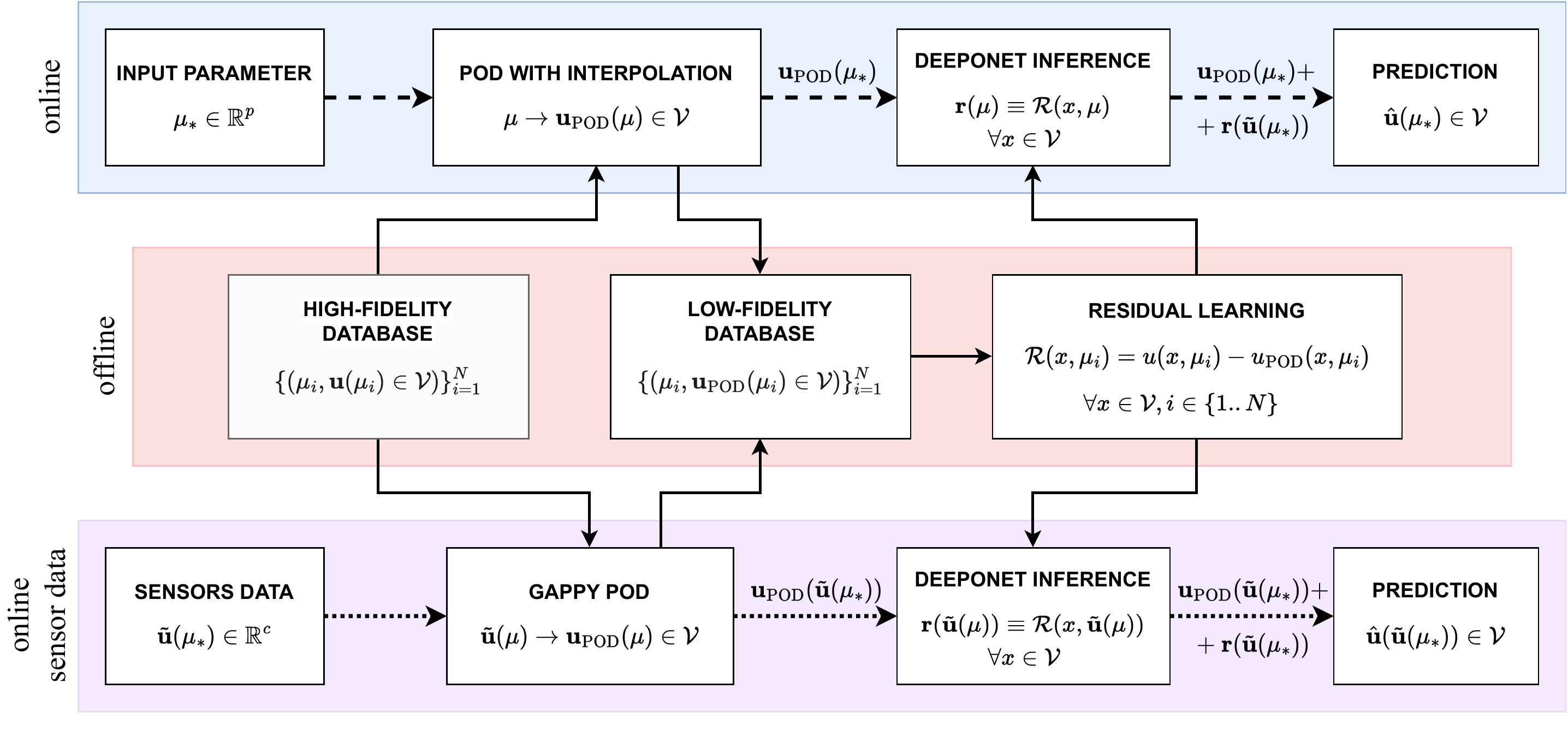}
  \caption{\RB{Scheme for the multi-fidelity POD framework. The arrows
      indicate the relationship between the different methods. In
      the blue box the dashed arrows indicate the online phase
      when the input parameter is provided. In the purple box the dotted arrows
      indicate the information flow if only sensor data are provided. The central
      red frame emphasizes the computationally expensive offline
      phase.}} \label{fig:scheme}
\end{figure}

For the prediction of solutions for new parameters, the non-intrusive POD model and the DeepONet are finally queried, as sketched in Figure~\ref{fig:scheme}. 
POD returns the low-fidelity (linear) approximation by providing the test parameter or sensor data, while the neural network returns the nonlinear residual. In some sense, this pipeline aims to exploit the advantages of the consolidated POD model, but at the same time improves it by adding a nonlinear term. So it can be also seen as a closure model. 

\subsection{Proper Orthogonal Decomposition for low-fidelity modeling}
\label{sec:pod}
POD is a consolidated technique widely used for model order
reduction. In this section, we briefly introduce how to compute the POD
modes, and we devote~\cref{sec:cpod} to present the gappy POD variant
in detail. 

The method consists of the computation of the optimal reduced basis to
represent the parametric solution manifold through a linear
projection. Let $u_i \in \mathbb R^n$ be the discrete solution corresponding
to the $i$-th parameter, and $U = [u_1, \dots, u_N] \in \mathbb{R}^{n
  \times N}$ be the snapshots matrix, whose columns are the solution
vectors. We want to find a linear approximation such that: 
\begin{equation}
u_i  \approx \sum_{k=1}^{r} a_i^k \psi_k, \qquad \text{ for } r\ll n, \text{ and for } i=1, \dots , N,
\end{equation}
where $\psi_i \in \mathbb{R}^n$ are the vectors comprising the reduced
basis, the so-called \emph{modes}, and $a_i := \begin{pmatrix}a_i^1,
  a_i^2, \dots, a_i^r\end{pmatrix} \in \mathbb R^r$ are the coordinates of
the corresponding solution at the reduced level, called modal
coefficients or latent variables. These reduced variables are obtained
by a projection of the solution snapshots onto the modes. \\
The POD modes can be obtained from the matrix $U$ in different ways:
by computing its singular value decomposition (SVD), or by decomposing
its correlation matrix~\cite{volkwein2013proper}. Moreover, all the
modes have a corresponding singular value, which represents their energetic
contribution. By arranging these
modes in decreasing order (with respect to the singular values), we
can express the original system with a hierarchical basis, from which
we can discard the less meaningful modes. \RA{The energy criterion based
on the singluar values decay reads as
\begin{equation}
\frac{\sum_{j=1}^r \sigma_j}{\sum_{j=1}^N \sigma_j} > \epsilon,
\end{equation}
where $\sigma_j$ is the $j$-th singular value, and $\epsilon$ is a
tolerance, usally set $\geq \num{0.99}$.} In other words, by providing
some samples of the solution manifold, POD is able to detect correlations between the data
and reduce the dimensionality of these discrete solutions. This
the approach becomes a fundamental tool for solving parametric partial
differential equations (PDEs) in a many-query context, mainly due to
the high-dimensional discrete spaces involved. 

The POD space can be exploited in a Galerkin framework, by projecting
the differential operators, or in a data-driven fashion by
coupling it with an interpolation (or regression) technique. In this
case, the database of reduced snapshots \RA{$\{\mu_i, a_i\}_{i=1}^N$} is
used as input to build the mapping $\mathcal I: P \to \mathbb R^r$ such that
$\mathcal I(\mu_i) = a_i$ for $i= 1, \dots, N$, which is
used for interpolating the modal coefficients for any new
parameter. Depending on the chosen regression technique, the equality
could not hold in principle, and we have $\mathcal I (\mu_i) \approx a_i$. Finally, exploiting such a mapping, we have the possibility to query for the modal coefficients at any test parameter belonging to the space $P$ and finally exploit the POD modes to map back the approximated solution in the original high-dimensional space.

\subsubsection{Gappy POD for sensors data}
\label{sec:cpod}
The main assumption for using gappy POD is to have access to only some sensor data. These sensors are placed at specific locations, given by the projection
matrix, or point measurement matrix, $C \in \mathbb{R}^{c \times n}$, with $c \ll n$, which contains
$1$ at measurements location and $0$ elsewhere. Using the
canonical basis vectors of $\mathbb{R}^n$ it takes the following form
\begin{equation}
C =
\begin{bmatrix}
e_{\gamma_1} & e_{\gamma_2} & \cdots & e_{\gamma_c} 
\end{bmatrix}^T,
\end{equation}
for some indices $\gamma_i \in [1, \dots, n]$, with $i \in [1, \dots, c]$.
The measurements $\tilde{u}_* \in \mathbb{R}^c$ of a generic full state vector
$u_* \in \mathbb{R}^n$ are thus given by
\begin{equation}
\tilde{u}_* = C u_*.
\end{equation}

If we now consider a parametric framework we can collect the parameter--solution
snapshot pairs $\{\mu_i, u_i \}_{i=1}^N$, where $\mu_i \in P \subset
\mathbb{R}^p$, and $u_i \in \mathbb{R}^n$ is the corresponding full
state. We arrange the snapshots by column in $U$ as
\begin{equation}
  U=\left[%
    \begin{array}{cccc}
      |   & |   & |   & |\\
      u_1 & u_2 & \dots & u_N\\
      |   & |   &  | & |
    \end{array}
  \right].%
\end{equation}
We take the $r$-rank SVD of the snapshots matrix $U$ and compute the
POD modes $\Psi_r$, so we can project the full states to their low-rank
representation $a \in \mathbb{R}^{r \times N}$:
\begin{equation}
U = \Psi \Sigma V^T \approx \Psi_r \Sigma_r V_r^T, \qquad
U \approx \Psi_r a.
\end{equation}
In the classical POD setting, where we deal with the full snapshots, we would just use the modal coefficients matrix $a$ to describe the solution manifold. For the gappy POD, instead, we have to consider the point measurement matrix. So, putting all together we have
\begin{equation}
\tilde{U} \approx (C \Psi_r) a,
\end{equation}
where $\tilde{U}$ is the matrix containing the sensors measurements $\{\tilde{u}_i \}_{i=1}^N$
arranged by columns, as done for the snapshots matrix.
For a generic snapshot $\tilde{u}_i$ we have:
\begin{equation}
\tilde{u}_i \approx C \sum_{k=1}^r  a_i^k \psi_k,
\end{equation}
where $\psi_k$ are the columns of $\Psi_r$, and $a_i^k$
are the modal coefficients, that is the $i$-th column of $a$.
A possible solution to find the modal coefficients is to minimize the residual in a least-squares
sense using the $L^2$ norm over the sensors locations which means
considering the following quantity~\cite{brunton2019data}
\begin{equation}
\int_{\text{supp}[\tilde{u}_i]} \left ( \tilde{u}_i - \sum_{k=1}^r  \tilde{a}_i^k \psi_k \right )^2.
\end{equation}

There are many ways in the literature to select the locations of the sensors:
optimal sensor locations that improve the condition number of $C
\Psi$~\cite{willcox2006unsteady, manohar2018data}, which are robust to
sensor noise, the sample maximal variance
positions~\cite{yildirim2009efficient}, or using information contained
in secant vectors between data points~\cite{otto2022inadequacy}, for
example. 
In this work, we are going to use the sparse sensor placement
optimization for reconstruction described in details
in~\cite{manohar2018data} and implemented in
PySensors~\cite{deSilva2021}.
The main idea is to find $C$ that minimizes the reconstruction error
using the modes $\Psi_r$ as in the following
\begin{equation}
C^* = \argmin_{C} \| U - \Psi_r (C \Psi_r)^\dagger \tilde{U} \|^2_2,
\end{equation}
where the symbol $\dagger$ stands for the Moore-Penrose pseudoinverse.

\subsection{DeepONet for residual learning}
\label{sec:multifidelity}
\begin{figure}[ht]
  \centering
\includegraphics[width=0.6\textwidth]{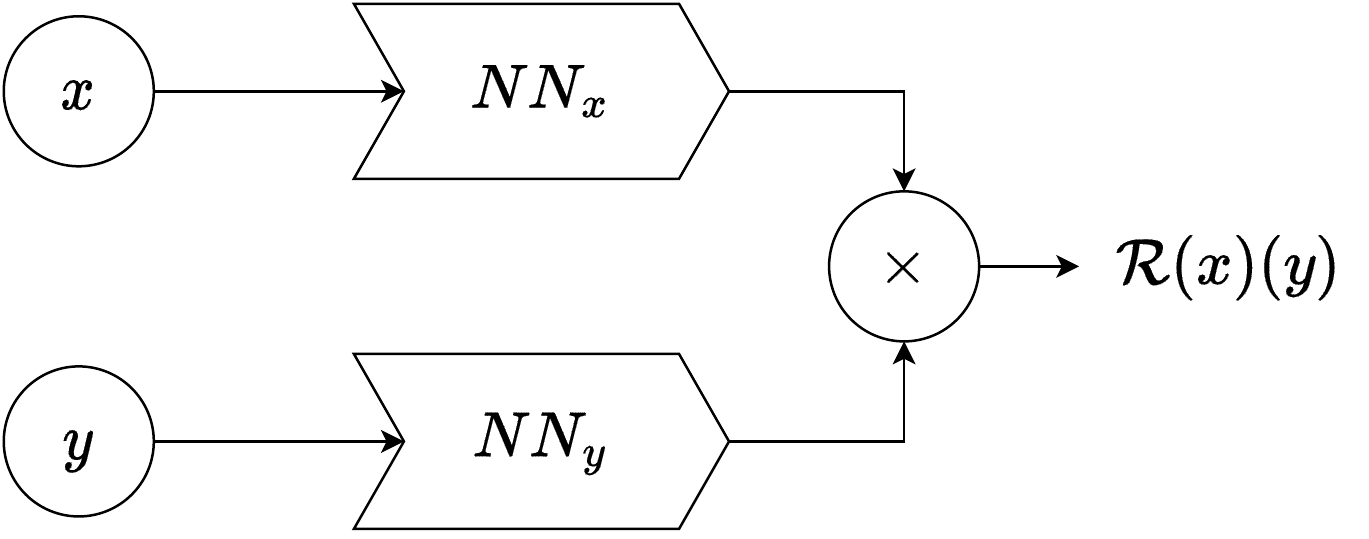}
\caption{The DeepONet scheme.\label{fig:deeponet}}
\end{figure}

DeepONet~\cite{lu2021learning} is a neural network architecture
able to learn nonlinear operators. Referring to the
original work for all the details, we emphasize its architecture composed
by two separate networks whose final outputs are multiplied to obtain
the final DeepONet outcome. The two networks, called
\emph{trunk} and \emph{branch}, can be any available architecture ---
e.g. convolutional network, graph network ---. In this work we 
consider feedforward networks (FFNs).
The networks are trained simultaneously during the learning loop: the
input is indeed divided into two independent components, $x \in \mathbb R^{N_x}$
and $y \in \mathbb R^{N_y}$, which feed the two networks $NN_x$ and $NN_y$,
respectively. The outputs $NN_x(x), NN_y(y) \in \mathbb R^{N_p}$ are finally
multiplied to approximate the operator $\mathcal G$:
\begin{equation}
\mathcal G(x)(y) \approx \sum_{i=1}^p [NN_x(x)]_i [NN_y(y)]_i.
\end{equation}
We underline that the choice of the two networks must satisfy the
dimensional constraint: they have to produce outputs with the same number
of components such that it is possible to compute their inner
product. The scheme in Figure~\ref{fig:deeponet} graphically
summarizes the structure of the DeepONet. 

\RB{In this work we adopt it to approximate the residual
$\mathcal{R}(x)(\mu) = u(x, \mu) - u_{\text{POD}} (x, \mu)$ in a multi-fidelity
approach}. We can
think at the mapping between the low-fidelity model (the POD/gappy POD)
and the high-fidelity model as a parametric operator $\mathcal
R(x)(\mu)$. This operator is numerically approximated by means of the
DeepONet, using as dataset the low- and high-fidelity databases
already computed. This architecture has demonstrated a great
capability in fighting overfitting issues~\cite{lu2021learning},
allowing to generalize the residual even with a limited set of
information.

\section{Numerical results}
\label{sec:results}
In this section we present the numerical results obtained by
applying the proposed numerical framework to a simple algebraic problem and to a
Navier-Stokes problem in a 2D domain. We are going to compare the proposed method with the POD model, the gappy POD model, and with the pure deep
learning approach by using DeepONet, aiming for a fair comparison with two state-of-the-art
techniques for (linear and nonlinear) data-driven modeling.
All the computations are
performed using PyTorch~\cite{pytorch} for the artificial neural networks,
EZyRB~\cite{demo2018ezyrb} for the POD with interpolation and gappy
POD calculations. To solve the Navier-Stokes equations with the finite element method we use FEniCS~\cite{logg2012automated}.

\subsection{Algebraic parametric problem}

The first test case is a simple benchmark problem inspired by~\cite{benamara2016multi}.
The high-fidelity parametric function $f^H : \Omega \times P \to \mathbb R$ is defined as
\begin{equation}
f^H(x; \mu) := \frac{1}{2} (\mu_1 x - 2)^2 \sin(12 x - 4) + \sin(\mu_2 \cos(5 x)),
\end{equation}
where $x \in \Omega = [0, 1] \subset \mathbb R$, and $\mu = (\mu_1, \mu_2) \in
P = [2, 15] \times [3, 20] \subset \mathbb R^2$. 
The first step is to compute the function value in some points in order to
build the high-fidelity database. We use different sampling strategies
for the spatial and parametric domain:
\begin{itemize}
\item we collect $n = 500$ equispaced samples $\{x_i^s\}_{i=1}^n$ in $\Omega$;
\item we collect $36$ samples using the latin hypercube
  sampling, plus $4$ additional samples at the corners of the domain, for
        a total of $N=40$ points $\{\mu_i^s\}_{i=1}^{N}$ in $P$.
\end{itemize}
We thus compose the snapshots matrix, varying the parametric
coordinates along the columns as follows:
\begin{equation}
\begin{bmatrix}
	f(x^s_1, \mu^s_1) & \dots & f(x^s_1, \mu^s_N)\\
	\vdots & \ddots & \vdots\\
	f(x^s_n, \mu^s_1) & \dots & f(x^s_n, \mu^s_N)\\
\end{bmatrix}
	\in \mathbb R^{n\times N}.
\end{equation}

Regarding the residual learning, we use the DeepONet model structured as
follows: the \emph{spatial} network (branch) is composed by $2$ inner
layers of $30$ neurons each, with the softplus activation function,
which is the smooth version of the Rectifier Linear Unit (ReLU)~\cite{glorot2011deep}; the
\emph{parametric} network (trunk) counts $2$ inner layers with $30$
neurons and the softplus function. The output layer has $30$
neurons for both networks, without applying any additional function
at this layer. The learning rate is equal to $\num{0.005}$, the
$L^2-$regularization factor is $\num{0.0001}$.

We propose a comparison between the MF approach, POD, and DeepONet in terms of accuracy on test parameters with a fixed input database of solutions.
We use different POD spaces in the comparison by selecting an
increasing energetic threshold for the modes selection, aiming to analyze the difference in the error
by varying the accuracy of the original POD model before getting
improved by MFDeepONet\footnote{For the remaining of this work, with \emph{MFDeepONet} we are going to
  refer to the proposed technique.}. \RA{We emphasize that no
  preprocessing or data centering is performed on the snapshots matrix,
resulting in the first mode representing a large amount of
energy. This corresponds to the minimal tolerance ($0.99$) in the
experiments below.}
Regarding the DeepONet architecture, we employ the one described above also to learn
the target function without the MF setting, such that the network
learns the actual unknown field instead of the residual. In this way,
we want to investigate the benefit of using the two methodologies (POD
and DeepONet) in a multi-fidelity fashion instead of only separately. 
We measure the relative error on an equispaced grid of $20\times 20$
parametric points.

\begin{figure}[htb]
  \centering
  \includegraphics[width=.995\textwidth,trim={3cm 1cm 2cm 1cm},clip]{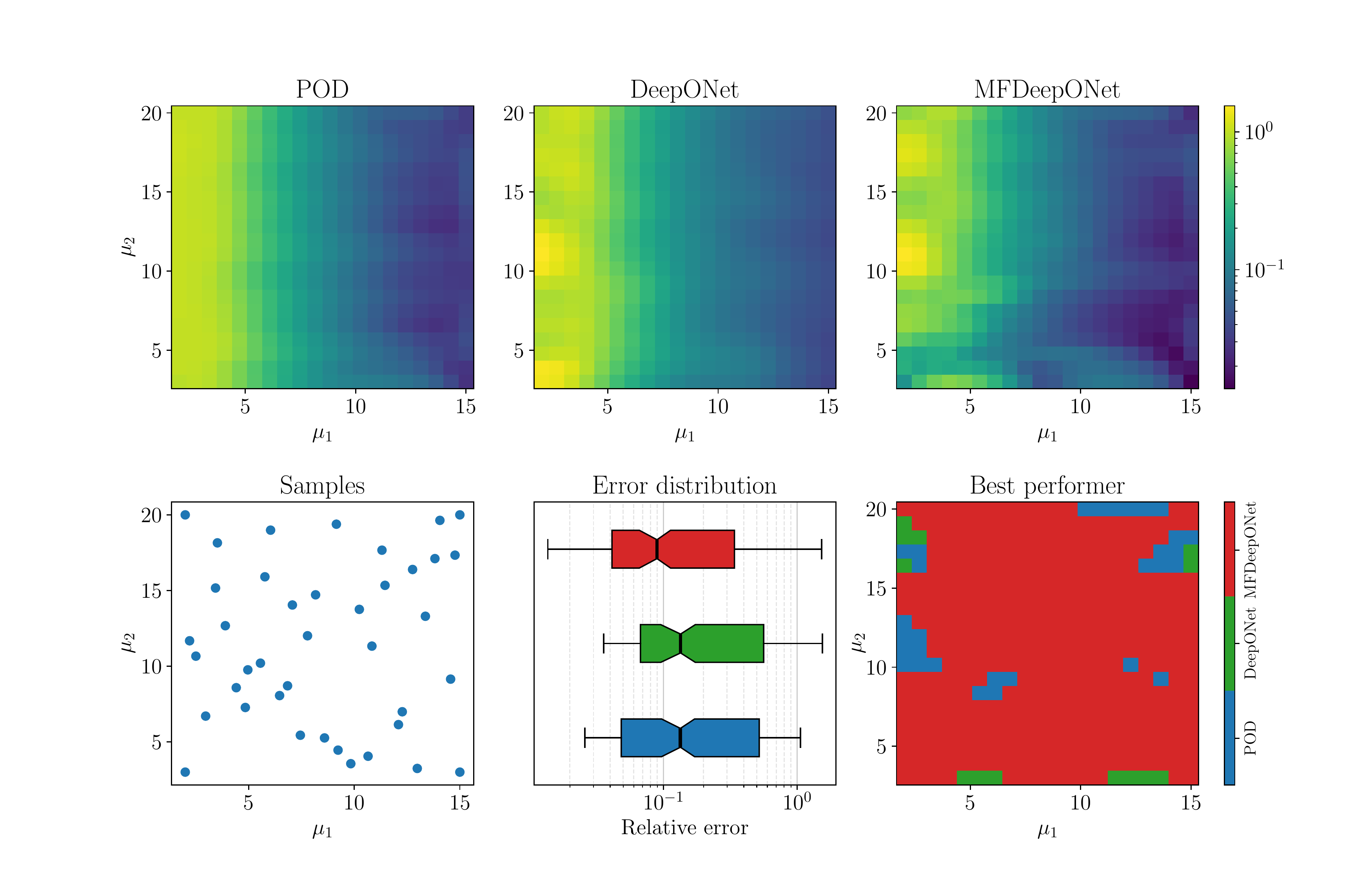}
  \caption{Comparison between POD ($0.99$ energy threshold), DeepONet, and multi-fidelity
	DeepONet. From top to bottom, we have the relative error in the
        parametric domain, the location of the high-fidelity samples in the
	parametric domain, the relative error distribution, and the best performers.}\label{fig:toy-comparison99}
\end{figure}

\paragraph{POD with energy threshold 0.99.}
For the POD model, we select an energy threshold $\epsilon =
\num{0.99}$ corresponding to $N = 1$ mode and radial basis function
(RBF) interpolation to approximate the map between the parameters and
the latent variables. The training for DeepONet and MFDeepONet lasts
$\num{10000}$ epochs. Figure~\ref{fig:toy-comparison99} shows a
quantitative comparison of the three investigated techniques,
presenting the relative error in the whole parametric domain, the
high-fidelity samples, and the error distribution. The last plot
(bottom right corner) graphically shows the technique which best
performs in all the tested parameters. 

In this experiment, the proposed methodology outperforms both POD and DeepONet. The relative error distribution suggests that mixing the techniques helps in terms of accuracy. 
\RB{Indeed, even if the error shows a greater variance, the MFDeepONet is able on average to achieve the best precision among the tested methods, resulting the better approach in almost all the parametric domain.}
We can also note that a direct correlation between the samples location and the error distribution is not visible, confirming the DeepONet capabilities in terms of generalization and making the proposed framework effective also during the testing phase.


\paragraph{POD with energy threshold 0.999.}
In this experiment, we replicate the previous settings with the exception of the new energy threshold for POD modes and a higher number of epochs for the machine learning models (DeepONet and MFDeepONet). Here we increase it to $\epsilon = \num{0.999}$ ($N = 6$ modes), addressing a more accurate original model, and balancing it with longer training.

\begin{figure}[htb]
  \centering
  \includegraphics[width=.995\textwidth,trim={3cm 1cm 2cm 1cm},clip]{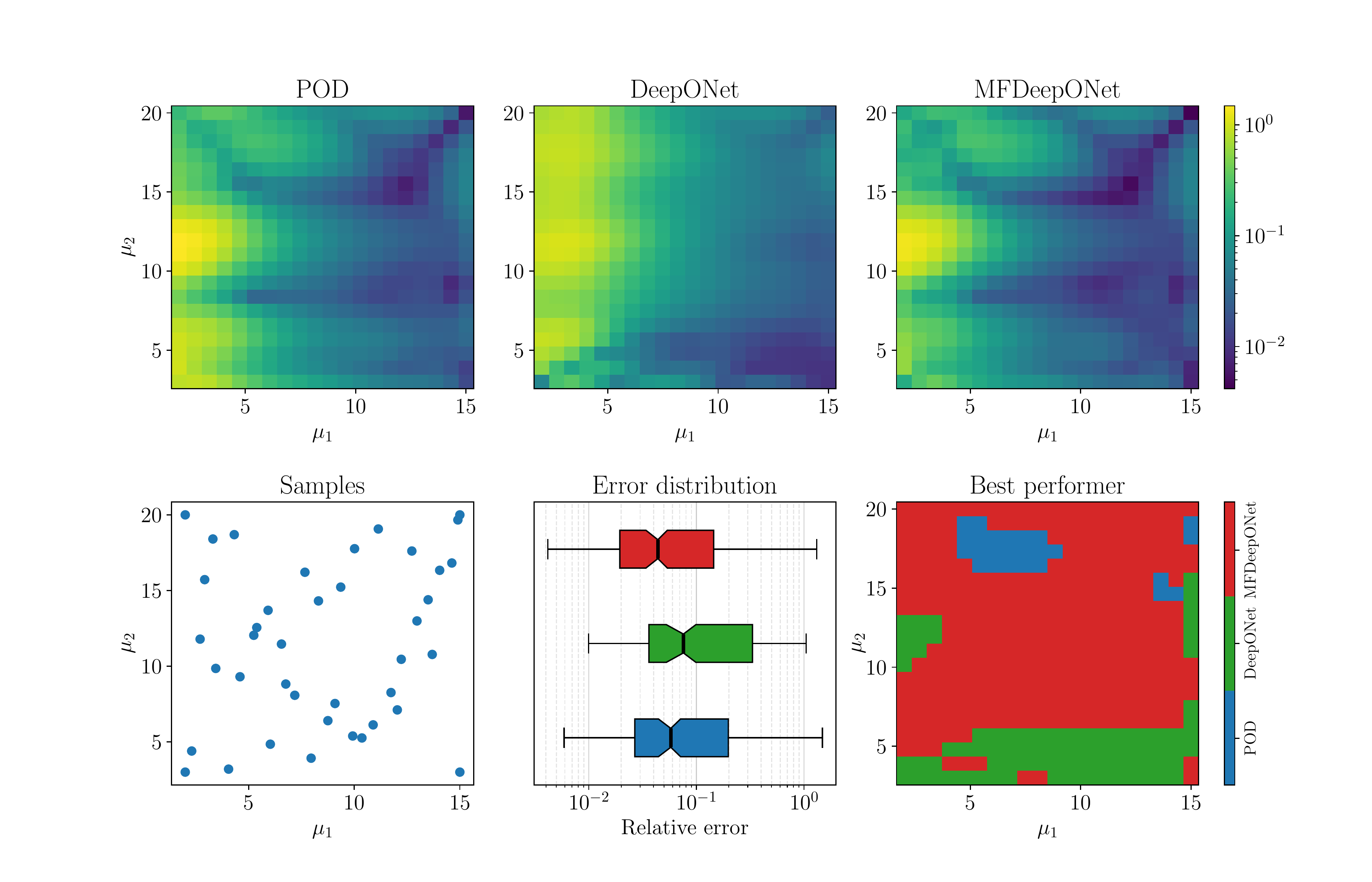}
  \caption{Comparison between POD ($0.999$ energy threshold), DeepONet, and multi-fidelity
	DeepONet. From top to bottom, we have the relative error in the
        parametric domain, the location of the high-fidelity samples in the
	parametric domain, the relative error distribution, and the best performers.}\label{fig:toy-comparison999}
\end{figure}

Figure~\ref{fig:toy-comparison999} illustrates the error obtained after a $\num{20000}$ epochs training. The results of the previous experiments are confirmed, even if with a lower overall benefit. The error distribution in the parametric space illustrates again how the MF enhancement combines the original methods: the regions of the parametric space where the methods work better are merged using MFDeepONet, resulting in a globally more accurate model.
However, using a more precise POD model (as low-fidelity) reduces the benefits of the MF approach, even with the higher number of epochs.


\paragraph{Gappy POD.}
Here we propose the same experiments as before, this time in a sensor
data scenario. Here we use $5$ sensor locations and a rank truncation
equal to $10$. 
The involved neural networks are trained in this case for $\num{50000}$ epochs.

\begin{figure}[htb]
  \centering
  \includegraphics[width=.995\textwidth,trim={3cm 1cm 2cm 1cm},clip]{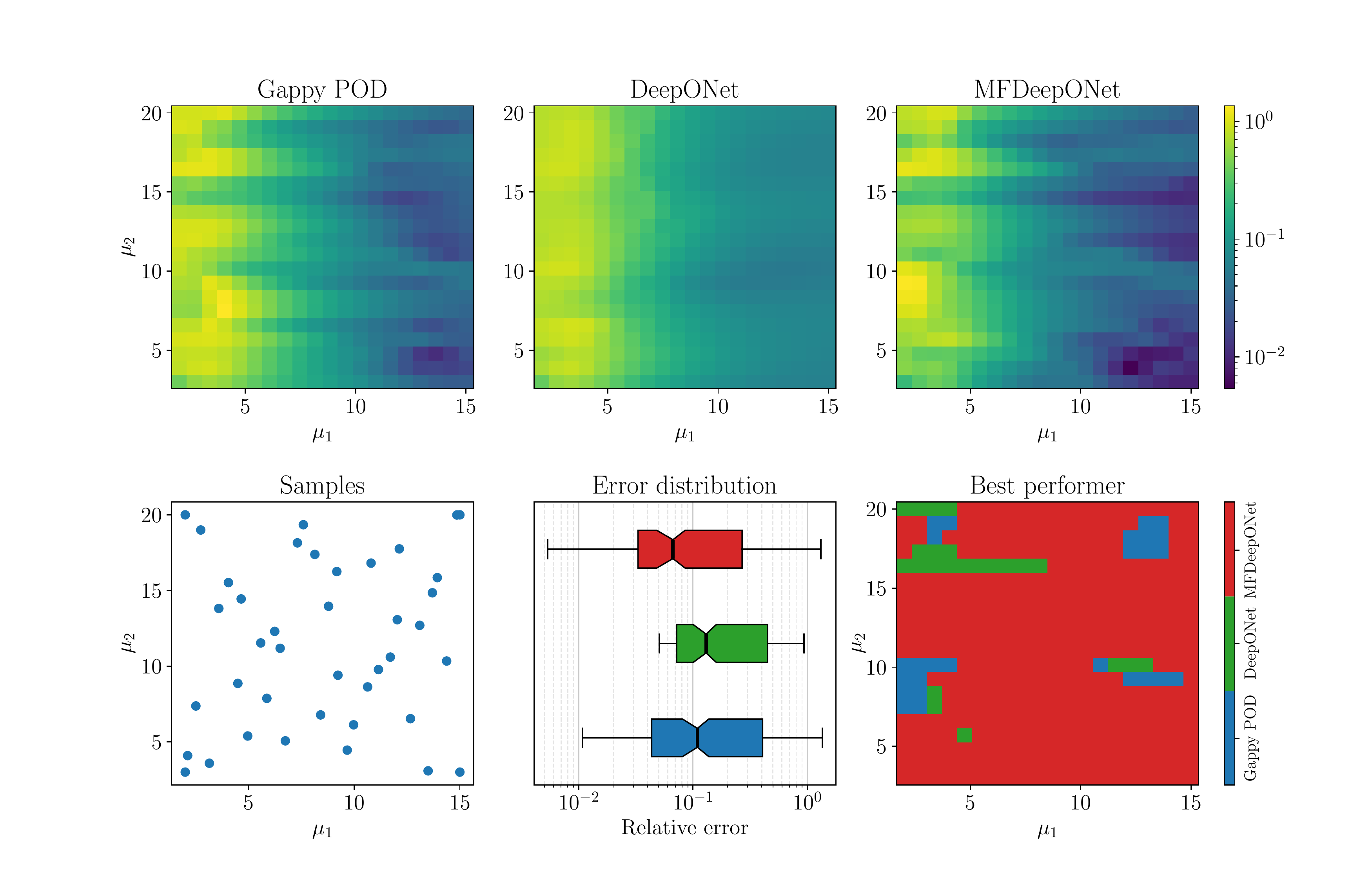}
  \caption{Comparison between gappy POD, DeepONet, and multi-fidelity
	DeepONet. From top to bottom, we have the relative error in the
        parametric domain, the location of the high-fidelity samples in the
	parametric domain, the relative error distribution, and the best performers.}\label{fig:toy-comparison}
\end{figure}

Figure~\ref{fig:toy-comparison} summarizes the
accuracy of the three tested methods, which are gappy POD, 
DeepONet, and MFDeepONet. The error distribution
demonstrates that the multi-fidelity approach performs statistically
better than the other methods. Looking at the competition between the
techniques, we can also note that the multi-fidelity approach reaches
the best accuracy in almost the whole parametric domain, even if at the
boundaries there is a precision decrease. Such an issue
could be mitigated by exploiting a better sampling strategy for the
high-fidelity data.

\begin{figure}[htb]
  \centering
    \includegraphics[width=.49\textwidth,trim={.5cm .5cm 0.5cm .5cm},clip]{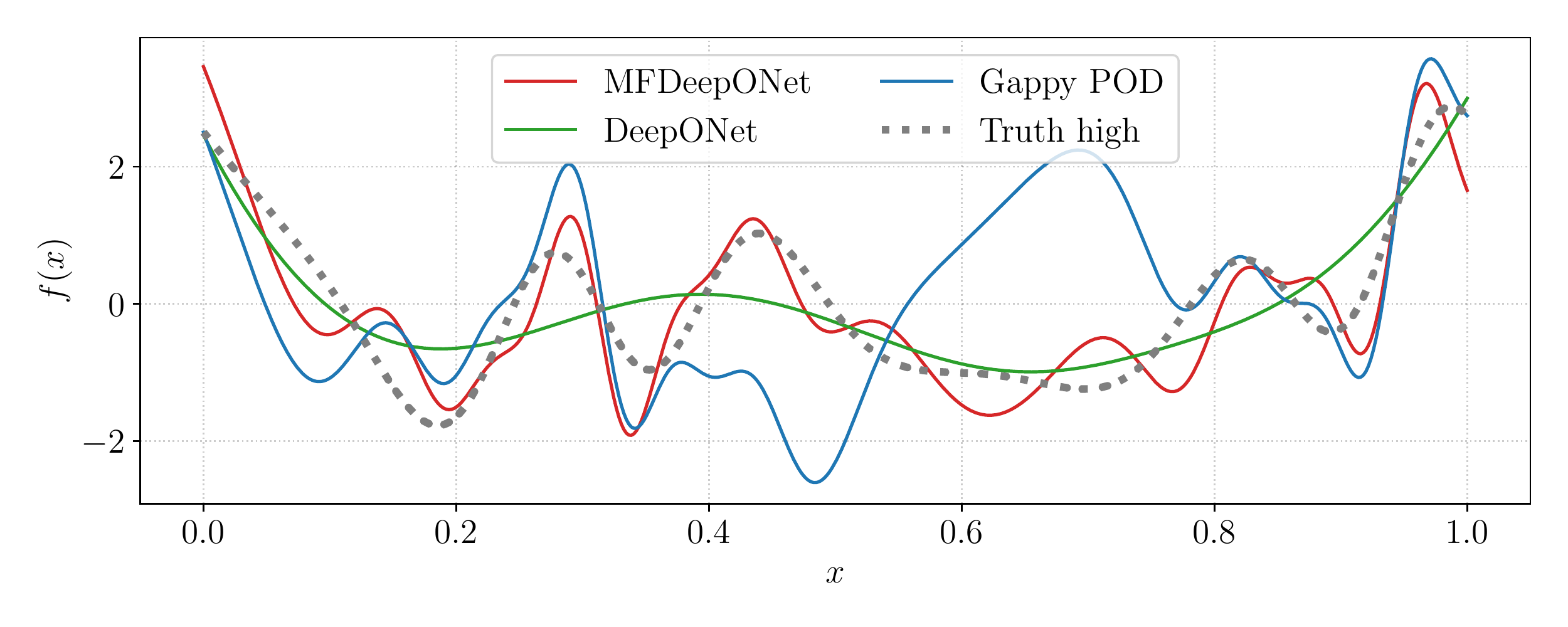}
    \includegraphics[width=.49\textwidth,trim={.5cm .5cm 0.5cm .5cm},clip]{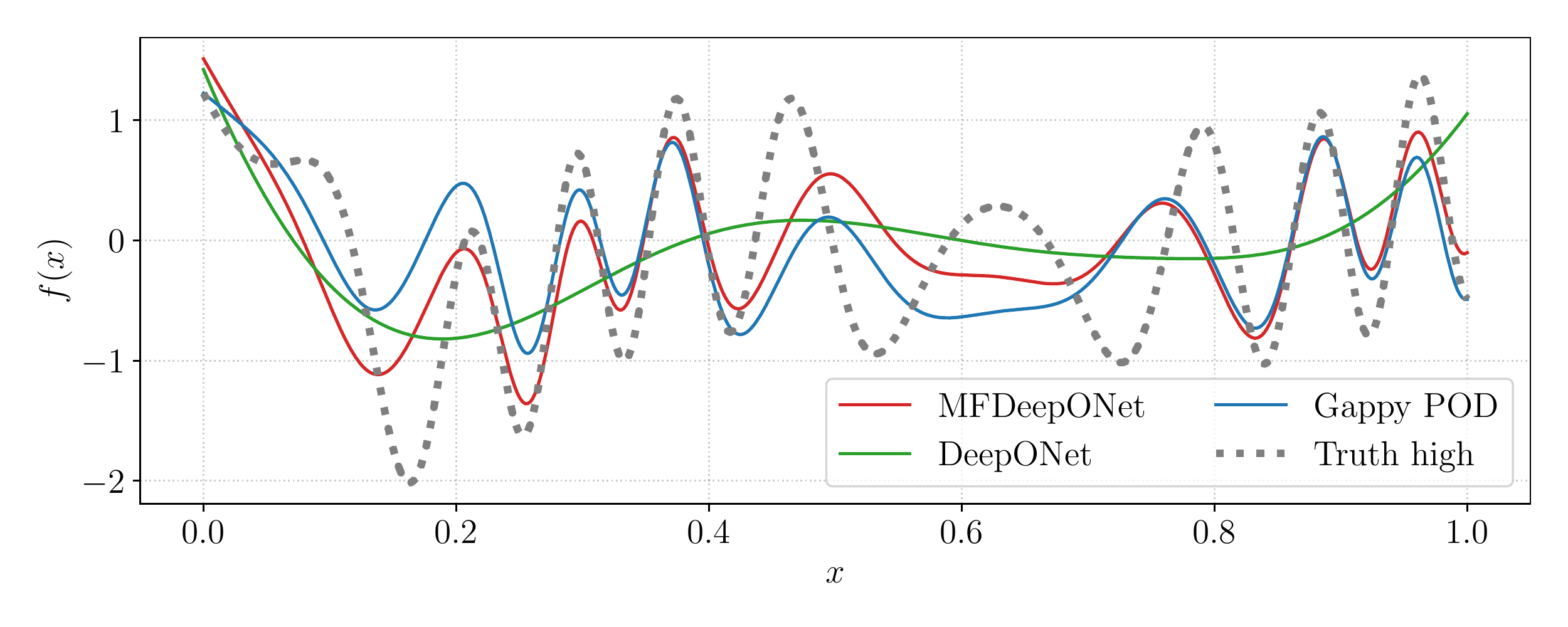}
    \includegraphics[width=.49\textwidth,trim={.5cm .5cm 0.5cm .5cm},clip]{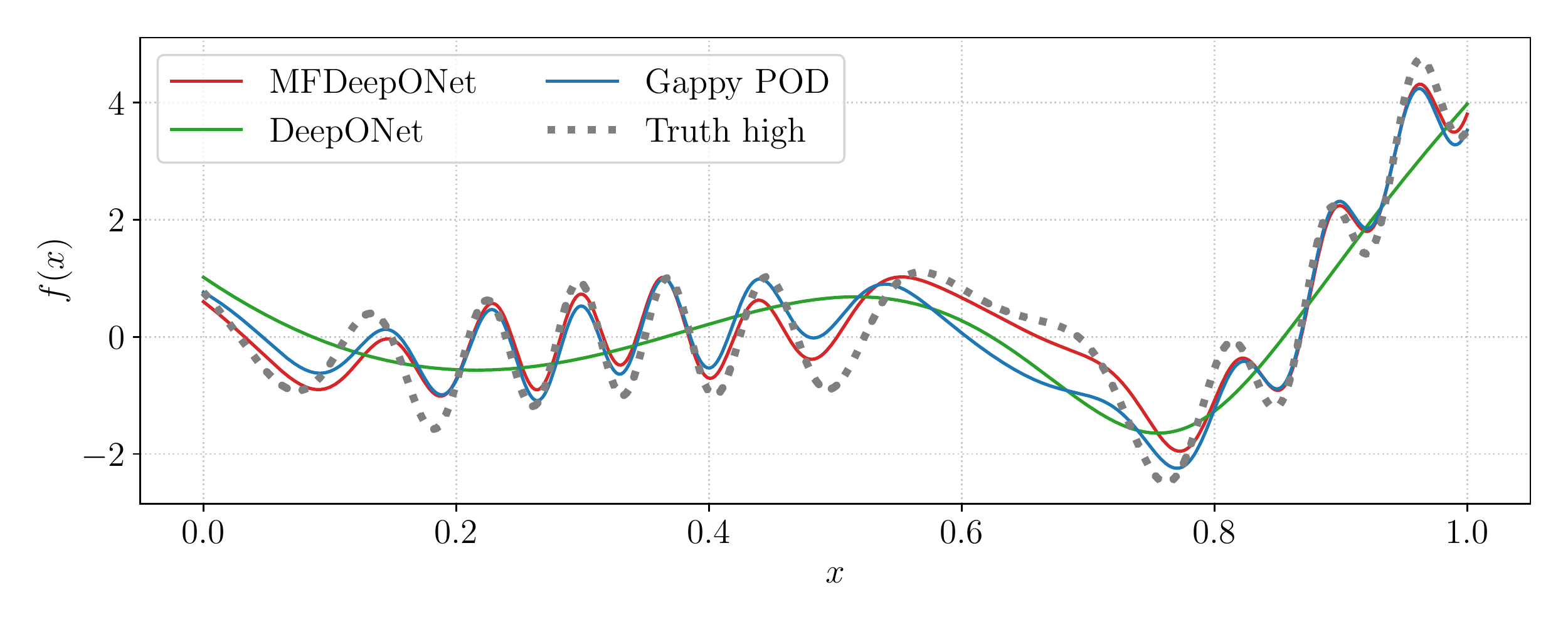}
    \includegraphics[width=.49\textwidth,trim={.5cm .5cm 0.5cm .5cm},clip]{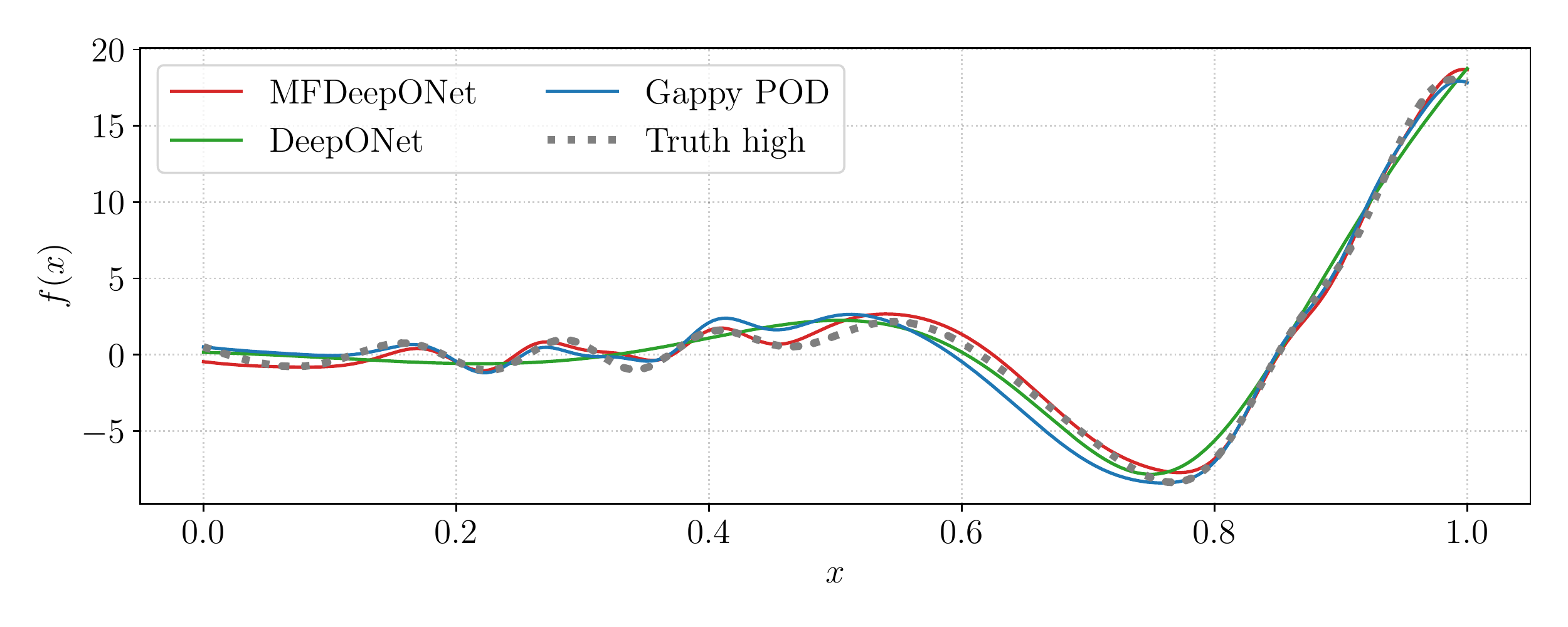}
	\caption{Examples of prediction using gappy POD, DeepONet and
	multi-fidelity DeepONet at different test parameters.
    \emph{Top left}: $\mu = [4, 8]$. 
    \emph{Top right}: $\mu = [3, 16]$.
    \emph{Bottom left}: $\mu = [5, 18]$.
    \emph{Bottom right}: $\mu = [8, 11]$.
    }\label{fig:toy518}
\end{figure}

The plots in Figure~\ref{fig:toy518} provide the comparison in
the spatial domain at four test parameters. The statistical
results are confirmed in these examples, with the multi-fidelity approach
that is able to predict most of the oscillations that the target function
exhibits, contrarily to the single-fidelity approaches.

\subsection{Navier Stokes problem}

In the second numerical experiment, we test the accuracy of the proposed method
for solving a parametric nonlinear PDE:
the incompressible Navier-Stokes equation on a 2D domain. The 
numerical setting is inspired by~\cite{BallarinManzoniQuarteroniRozza2015}.

\begin{figure}[htb]
  \centering
  \includegraphics[width=.8\textwidth]{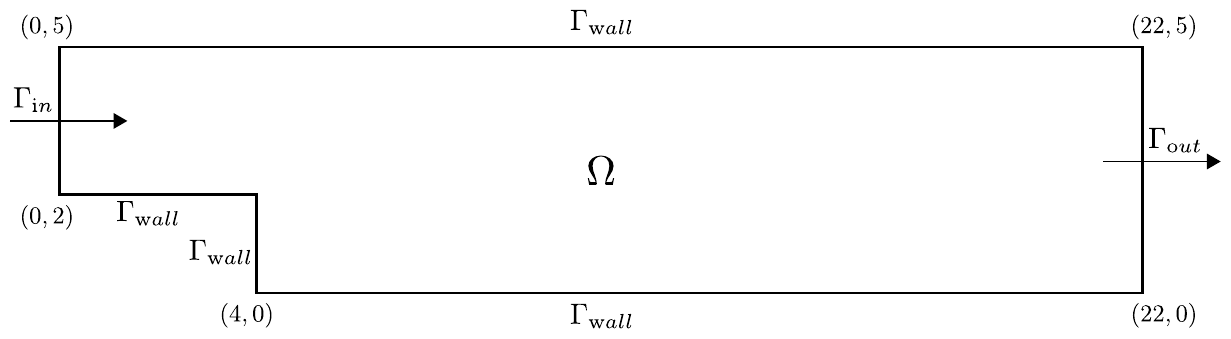}
  \caption{Domain description.}
  \label{fig:domain_ns}
\end{figure}

We define the parametric vector field $u: \Omega\times P \to \mathbb R^2$ and the parametric scalar field $p : \Omega\times P \to \mathbb R$ such that:
\begin{equation}
\left\{ 
\begin{array}{rll}
\nu \Delta \mathbf{u} + (\mathbf{u} \cdot \nabla)\mathbf{u} + \nabla p &= 0\quad\quad\quad&\text{in}\,\Omega,\\
\nabla \cdot \mathbf{u} &= 0\quad\quad\quad&\text{in}\,\Omega,\\
\mathbf{u} &=\mu\bigl\{\frac{1}{2.25}(x_1-2)(5-x_1),0\bigr\} \quad\quad\quad&\text{on}\,\Gamma_\text{in},\\
\mathbf{u} &= 0\quad\quad\quad&\text{on}\,\Gamma_\text{wall},\\
\nu \frac{\partial \mathbf{u}}{\partial \mathbf{n}} - p\mathbf{n} &= 0\quad\quad\quad&\text{on}\,\Gamma_\text{out},\\
\end{array}
\right.
\label{ns}
\end{equation}
where $x = (x_0, x_1) \in \Omega \subset \mathbb R^2$ and $\mu \in P = [1, 80]$.
The \emph{L}-shape spatial domain $\Omega$, together with the boundaries,
is sketched in Figure~\ref{fig:domain_ns}. For this test case, the
parametric solution is computed numerically by means of finite element
discretization. The spatial domain has been tessellated into $1639$
non-overlapping elements, and for stability we apply the Taylor-Hood
$P2-P1$ scheme. The high-fidelity dataset is composed of $20$ equispaced
parametric samples in $P$, arranged in the snapshots matrix $U \in
\mathbb R^{n\times N}$ with $N = 20$ and $n = 1639$.

The DeepONet structure for this problem is the following:
\begin{itemize}
\item the spatial network (branch) is composed by $3$ hidden
		layers of $50$ neurons each;
\item the parameter network (trunk) is composed by $3$ hidden
		layers of $20$ neurons each.
\end{itemize}
Also in this case, the last layer of the networks has the same number of
neurons, $20$. The activation function used in all the hidden layers
is the Parametric ReLU (PReLU)~\cite{he2015delving}, with the learning
rate equal to $\num{0.003}$ and the $L^2-$regularization factor equal to
$\num{0.001}$. The learning phase lasted $2.5 \times 10^4$ epochs. 
The accuracy of the MF approach is compared to the gappy POD and to
the standard DeepONet, with the same architecture (single-fidelity). 
The relative error is evaluated over $500$ testing points, randomly
sampled in the parametric space. 

\paragraph{POD with energy threshold 0.99.}
As before, we start with a relatively poor POD model, using $N=1$ mode selected by the energetic criterion. RBF is employed also here to approximate the solution manifold at the reduced level. The number of epochs is fixed at $\num{10000}$ for the deep learning training.

\begin{figure}[htb]
  \centering
  \includegraphics[width=.99\textwidth,trim={.5cm .5cm .5cm 1cm},clip]{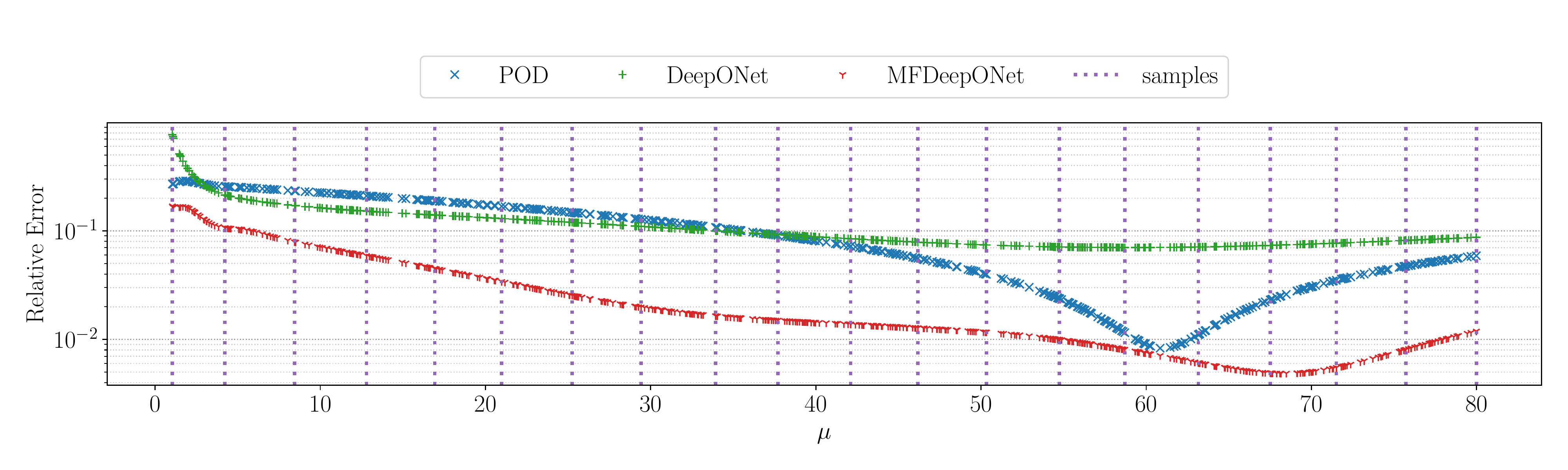}
  \caption{Comparison between POD (0.99 energy threshold), DeepONet, and multi-fidelity
	DeepONet in terms of relative error in the parametric domain. The vertical dotted lines indicate the location of the high-fidelity samples.}\label{fig:compare_ns99}
\end{figure}

Figure~\ref{fig:compare_ns99} shows the plot of the mean relative error over the spatial domain for all the test parameters, reporting also the location of the samples in the parameter space.
As for the previous experiment, the proposed technique is able to keep a
higher precision in the entire domain, without showing a visible
correlation between the location of the high-fidelity data and the error
trend, demonstrating its robustness in terms of possible overfitting.
Employing the DeepONet architecture to learn the residual (between the
POD and high-fidelity models) rather than the target function results in
a more efficient learning procedure, capable to \RA{ourperform the single-fidelity approaches in the entire parametric space here considered}.

\paragraph{POD with energy threshold 0.999.}
As for the previous test case, we repeat the same experiment with a more accurate POD model. Here we use $N=3$, raising the training time to $\num{20000}$ epochs.

\begin{figure}[htb]
  \centering
  \includegraphics[width=.99\textwidth,trim={.5cm .5cm .5cm 1cm},clip]{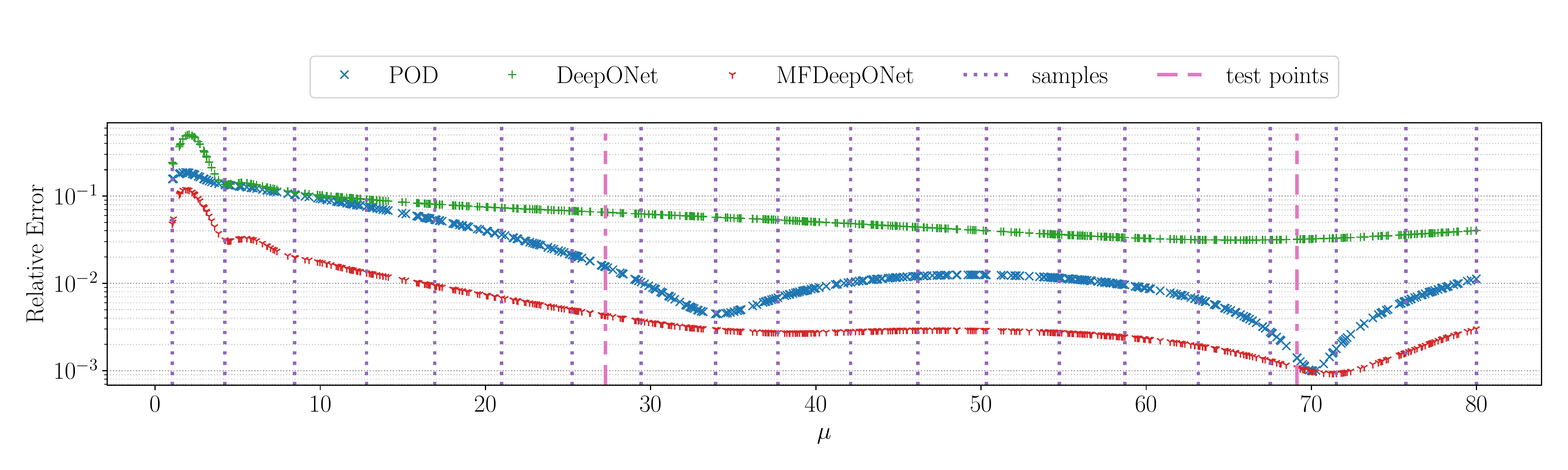}
	\caption{Comparison between POD (0.999 energy threshold), DeepONet, and multi-fidelity
	DeepONet in terms of relative error in the parametric domain. The vertical dotted lines indicate the location of the high-fidelity samples. The $2$ dashed vertical lines indicate the test parameters represented in Figures~\ref{fig:testns1}~and~\ref{fig:testns2}.}\label{fig:compare_ns999}
\end{figure}

The trend showed in the previous investigations is confirmed, as depicted in Figure~\ref{fig:compare_ns999}. The MFDeepONet method is able to produce a more accurate prediction in all the testing points, with no visible correlation with the training data.
For a fair comparison, we also investigated the predicted field in the only
point of the parametric domain where the MFDeepONet shows a slightly
higher error with respect to the POD model (whereas the standard
DeepONet performs poorly there). 

\begin{figure}[htb]
  \centering
  \includegraphics[width=.99\textwidth,trim={3cm .5cm 3cm 1cm},clip]{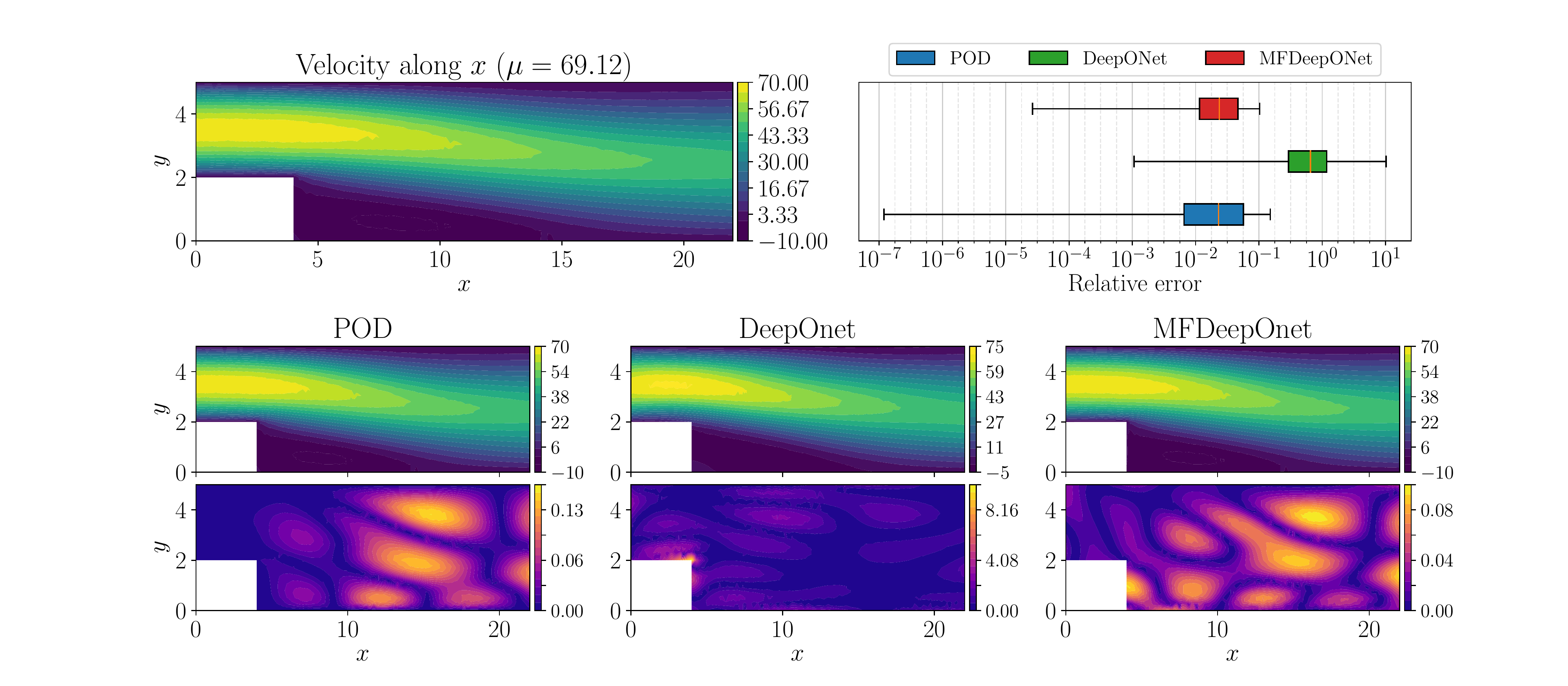}
	\caption{Representation over the spatial domain of the velocity (along $x$) in the Navier-Stokes testcase for $\mu = 69.12$. The approximation computed by POD (0.999 energy threshold), DeepONet, and multi-fidelity DeepONet is shown at the bottom together with the relative error. The distribution of the error is summarized in the box plot.}\label{fig:testns1}
\end{figure}

\begin{figure}[htb]
  \centering
  \includegraphics[width=.99\textwidth,trim={3cm .5cm 3cm 1cm},clip]{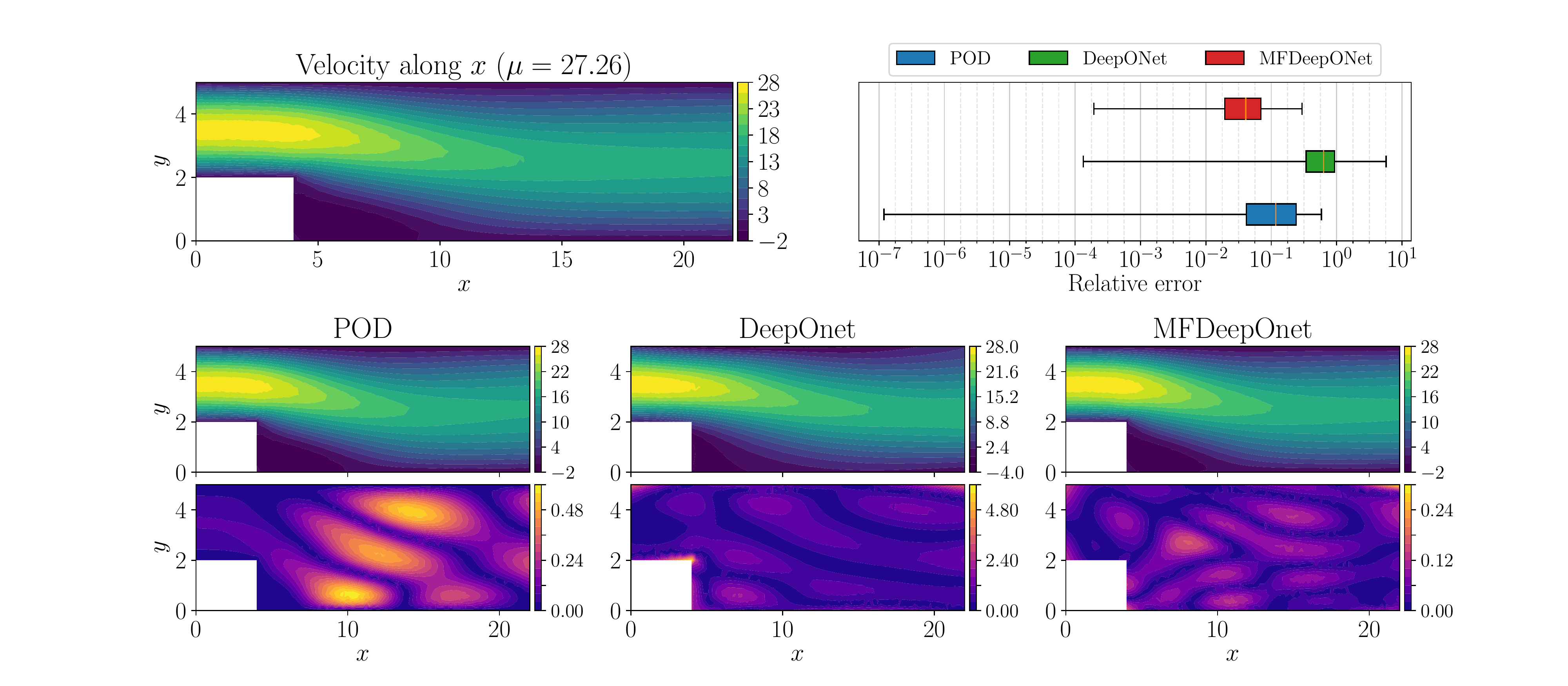}
	\caption{Representation over the spatial domain of the velocity (along $x$) in the Navier-Stokes testcase for $\mu = 39.95$. The approximation computed by POD (0.999 energy threshold), \RA{DeepONet, and multi-fidelity} DeepONet is shown at the bottom together with the relative error. The distribution of the error is summarized in the box plot.}\label{fig:testns2}
\end{figure}

Figure~\ref{fig:testns1} shows the $x$-component of the velocity field for the parameter $\mu = 69.12$ obtained by the three methods, with a statistical summary of the relative error. The MF approach shows here a smaller spatial variance, even if on average performs equally to the POD model.
Looking instead at a different parametric coordinate
(Figure~\ref{fig:testns2}), the benefits of the proposed approach
become clear. The considerations regarding the variance of the error
are still valid, but the solution for $\mu = 39.95$ shows \RA{a
  remarkable improvement in the accuracy over the testing points}.

\paragraph{Gappy POD.}
The last numerical experiment focuses on the Navier-Stokes model, for which
sensor data are used by the gappy POD for the low-fidelity
approximation. Here we use $7$ sensor locations and a rank truncation
equal to $8$.
We trained the DeepONet and the MFDeepONet for $\num{50000}$ epochs.

\begin{figure}[htb]
  \centering
  \includegraphics[width=.99\textwidth,trim={.5cm .5cm .5cm 1cm},clip]{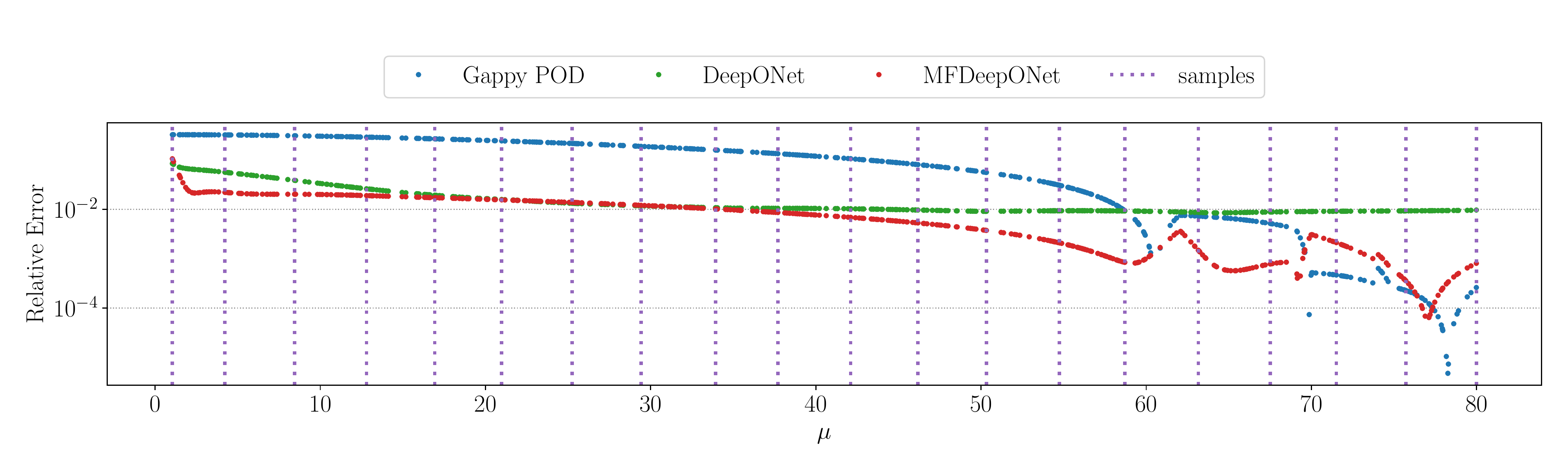}
	\caption{Comparison between gappy POD, DeepONet, and multi-fidelity
	DeepONet in terms of relative error in the parametric domain. The vertical dotted lines indicate the location of the high-fidelity samples.}\label{fig:compare_ns}
\end{figure}

Figure~\ref{fig:compare_ns} reports the relative test error measured in all the test points. In this case, the standard DeepONet, is able to outperform the POD model in a large region of the parametric domain, with a relative error that remains close to $\num{0.01}$. Gappy POD is able to reach the best precision in a few test points, but also here the MF approach is the best compromise in terms of global accuracy, even if it is actually less precise than the POD model for high parameter values ($\mu > 70$).

\subsection{Summary discussion}
This section is devoted to a summary discussion of the results obtained in the numerical investigations. For a fair comparison, we computed the mean relative test error\footnote{We recall the test error is computed over a $20\times20$ regular grid for the algebraic problem, and at $460$ random sample for the Navier-Stokes problem.} for each method, reporting the accuracy for different neural networks training times. In addition to the previous tests, we show in Table~\ref{tab:summary} the results obtained by employing a POD space whose modes are selected with an energetic threshold of $\epsilon = \num{0.9999}$. 
The error charts for the missing cases, as well as some graphical representations of the parametric solutions, are reported in Appendix~\ref{sec:appendix}.
The latter experiment aims to analyze the final accuracy when the low-fidelity POD is even more precise: the Mf approach is able to reach the best mean relative error, but its effectiveness is marginal, confirming the trend already defined in the previous tests. The combination of the POD model and DeepONet in the cascade fashion is able to reach the best accuracy in almost all the cases, but its improvement becomes marginal when the POD has good accuracy. Learning the residual however does not seem to affect the final outcome in a pejorative way, provided that the DeepONet is trained for a proper number of epochs. This is for sure a critical issue inherited by deep learning in general: we can indeed see that a longer training step does not always ensure better accuracy, producing instead over-fitting. On the practical side, the optimal settings of the network --- e.g. training epochs, number of layers, type of activation function --- need to be calibrated with a trial and error procedure or using more sophisticated approaches such as grid search. This calibration is out of the scope of this investigation where we want to formalize the novel framework, but surely sensitivity analysis regarding the hyper-parameters will be explored in future works.
The generalization of the DeepONet, assisted also by the $L_2$-regularization imposed during the optimization, is able to improve accuracy over the entire parametric space, without showing a visible correlation between the location of the high-fidelity snapshots and the relative error spatial distribution.

To conclude, we highlight that the numerical experiments demonstrate a great improvement when the original POD model lacks accuracy, resulting in a great tool to treat problems where POD is not able to capture all the fluid characteristics, due to the complexity of the mathematical model or to the limited number of high-fidelity snapshots. 

\begin{table}[]
    \centering
    \caption{The mean relative error computed in all the experiments. In bold the best results for each row.}\label{tab:summary}
    \begin{tabular}{llccccccc}
         \toprule
         & & \multirow{2}{*}{POD} & \multicolumn{3}{c}{DeepONet} &\multicolumn{3}{c}{MFDeepONet}\\
         & & & 10k & 20k & 50k & 10k& 20k& 50k\\
         \midrule
         \multirow{3}{*}{testcase \#1} 
         & POD rank = 0.99   & 0.324 & 0.270 & 0.265 & \textbf{0.217} & 0.293 & 0.247 & 0.308\\
         & POD rank = 0.999  & 0.203 & 0.270 & 0.265 & 0.217 & 0.196 & 0.193 & \textbf{0.127}\\
         & POD rank = 0.9999 & 0.098 & 0.270 & 0.265 & 0.217 & \textbf{0.093} & 0.104 & 0.178\\
         \midrule
         \multirow{3}{*}{testcase \#2} 
         & POD rank = 0.99   & 0.105 & 0.116 & 0.068 & 0.072 & 0.030 & \textbf{0.022} & 0.030\\
         & POD rank = 0.999  & 0.033 & 0.116 & 0.068 & 0.072 & 0.025 & 0.023 & \textbf{0.009}\\
         & POD rank = 0.9999 & 0.011 & 0.116 & 0.068 & 0.072 & 0.011 & \textbf{0.009} & \textbf{0.009} \\
         \bottomrule
    \end{tabular}\\[0.4cm]
    \begin{tabular}{llcccccc}
         \toprule
         & \multirow{2}{*}{Gappy POD} & \multicolumn{3}{c}{DeepONet} &\multicolumn{3}{c}{MFDeepONet}\\
         & & 10k & 20k & 50k & 10k& 20k& 50k\\
         \midrule
         testcase \#1 & 0.260 & 0.419 & 0.380 & 0.278 & 0.218 & 0.217 & \textbf{0.197} \\
         testcase \#2 & 0.135 & 0.035 & 0.033 & 0.017 & 0.029 & 0.010 & \textbf{0.009} \\
         \bottomrule
    \end{tabular}
    \label{tab:my_label}
\end{table}

\section{Conclusions and future perspectives}
\label{sec:conclusions}
In this work, we introduced a novel approach to enhance POD-based 
reduced order models thanks to a residual learning procedure by DeepONet.
It operates by building from a limited set of data an initial low-fidelity approximation exploiting
established reduced order modeling techniques. Then it learns the difference
between this low-fidelity representation and the original model through
the artificial neural networks, that will be inferred to predict the solution at unseen parameters. We
emphasize that such an enhancement neither needs any additional
evaluation of the original model nor the knowledge of the
high-fidelity model, resulting in a generic data-driven improvement
at a fixed computational budget.
This framework has demonstrated its effectiveness in two different
testcases: a univariate parametric function and a Navier-Stokes problem
on a $2$-dimensional domain, showing a higher precision in both 
experiments with respect to the use of single-fidelities. We highlight that in these experiments the number of considered POD modes is voluntarily kept small, simulating a POD model with poor accuracy.

    The present work illustrates the pipeline for POD and gappy POD for the
construction of the low-fidelity model and the DeepONet architecture for
residual learning. Due to its modularity, the framework is general, admitting in principle to replace the low-fidelity models with different
ones. Possible future extensions should investigate adaptive samplings
and sensor placement exploiting the proposed numerical framework.

\section*{Acknowledgements}
This work was partially supported by an industrial Ph.D. grant sponsored by
Fincantieri S.p.A. (IRONTH Project), by the MIT-FVG project ``Multi-disciplinary Ship Design by Reduced Order Models and Machine Learning'', and partially funded by European
Union Funding for Research and Innovation --- Horizon 2020 Program --- in the framework of European Research Council Executive Agency: H2020 ERC CoG 2015
AROMA-CFD project 681447 ``Advanced Reduced Order Methods with Applications in
Computational Fluid Dynamics'' P.I. Professor Gianluigi Rozza.

\section*{Appendix}
\label{sec:appendix}
\begin{figure}[H]
    \centering
    \includegraphics[width=.49\textwidth,trim={.5cm .5cm 0.5cm .5cm},clip]{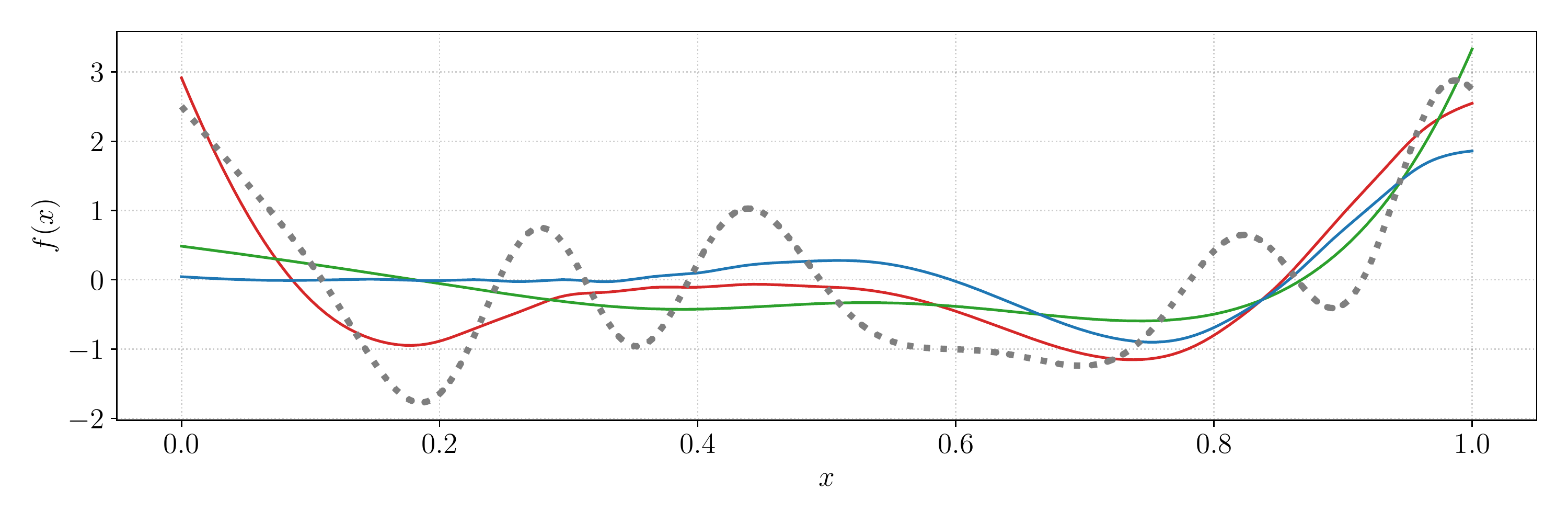}
    \includegraphics[width=.49\textwidth,trim={.5cm .5cm 0.5cm .5cm},clip]{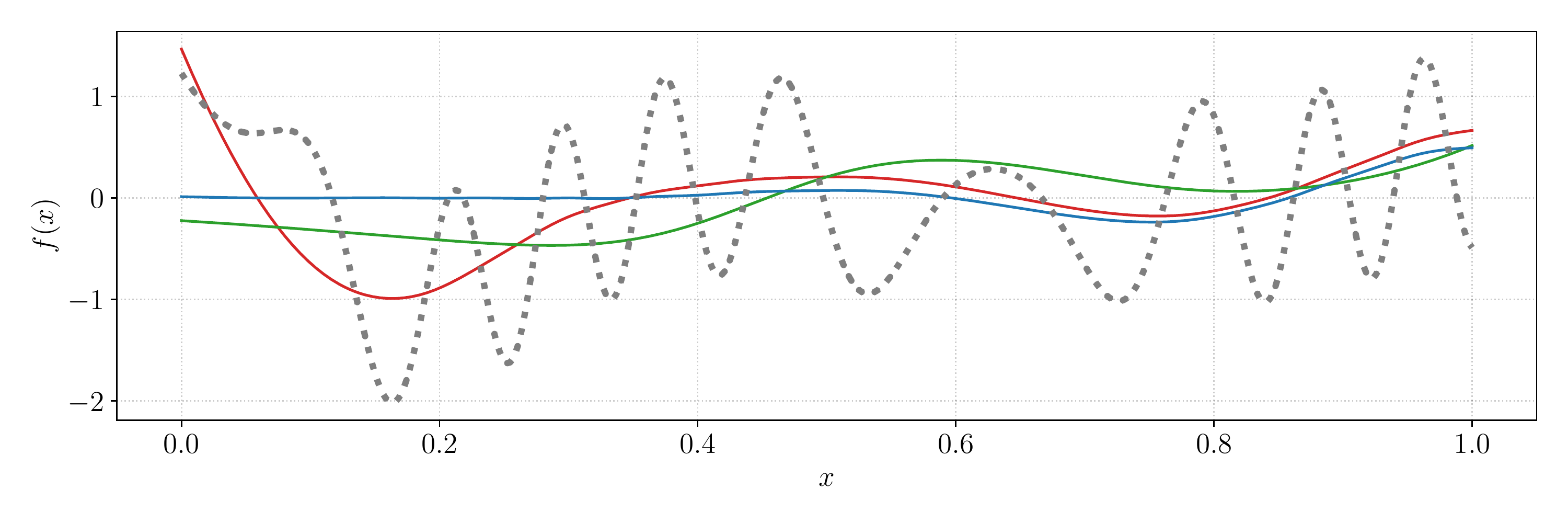}
    \includegraphics[width=.49\textwidth,trim={.5cm .5cm 0.5cm .5cm},clip]{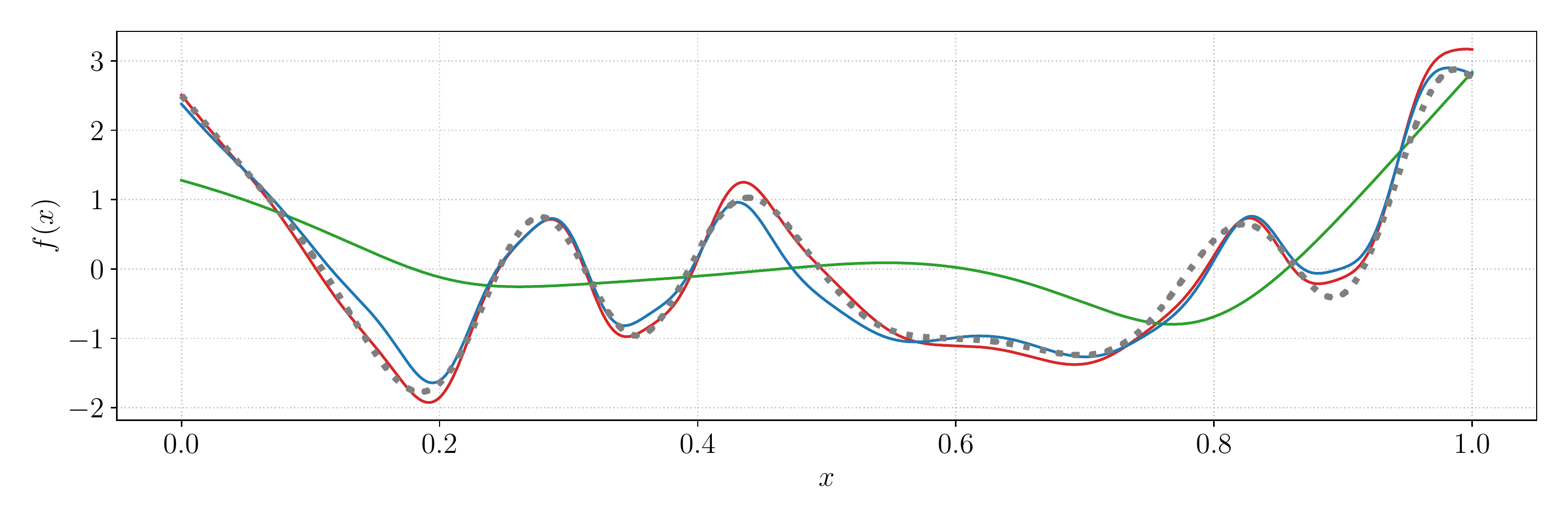}
    \includegraphics[width=.49\textwidth,trim={.5cm .5cm 0.5cm .5cm},clip]{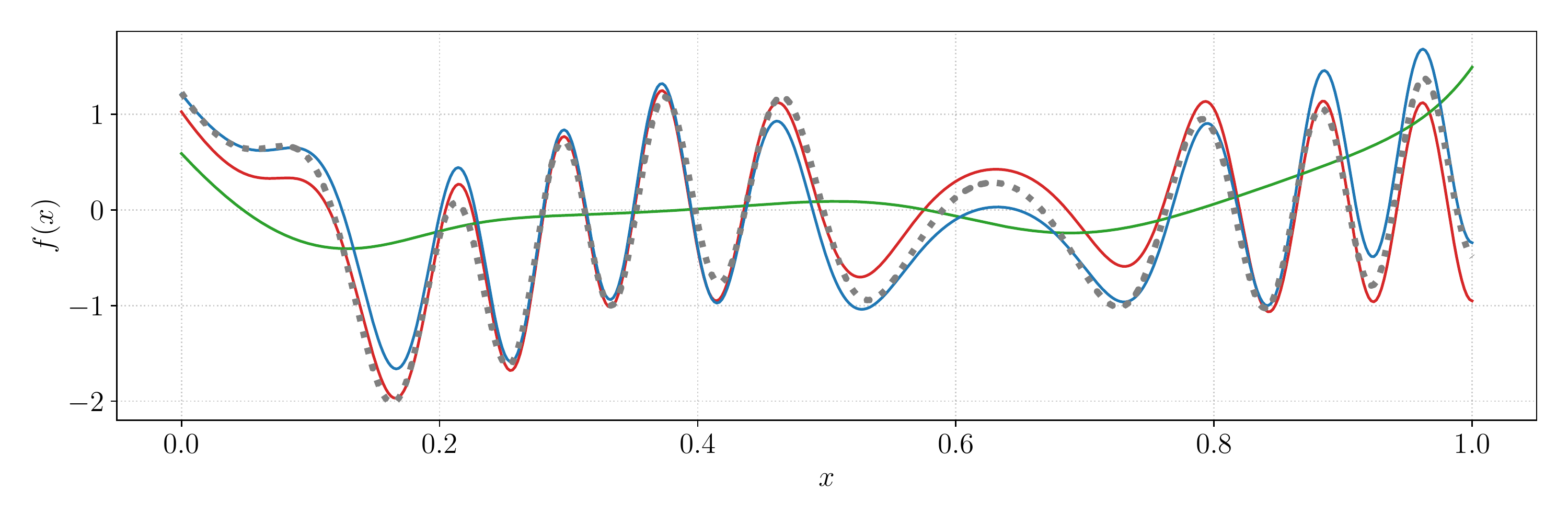}
    \includegraphics[width=.49\textwidth,trim={.5cm .5cm 0.5cm .5cm},clip]{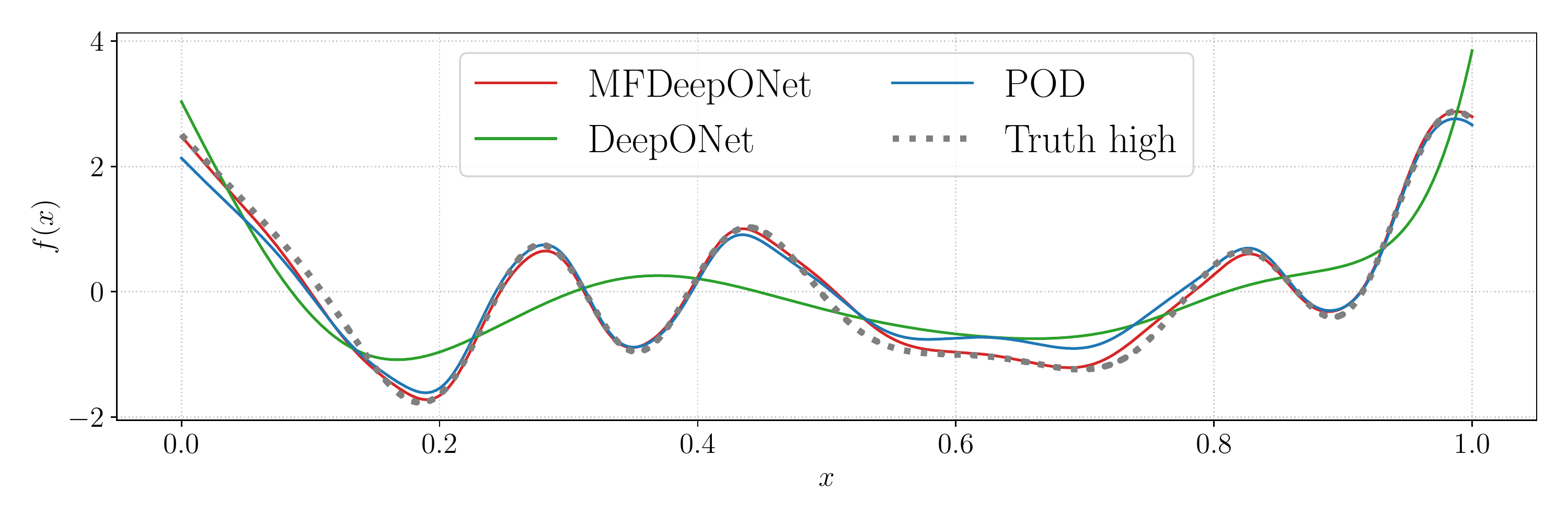}
    \includegraphics[width=.49\textwidth,trim={.5cm .5cm 0.5cm .5cm},clip]{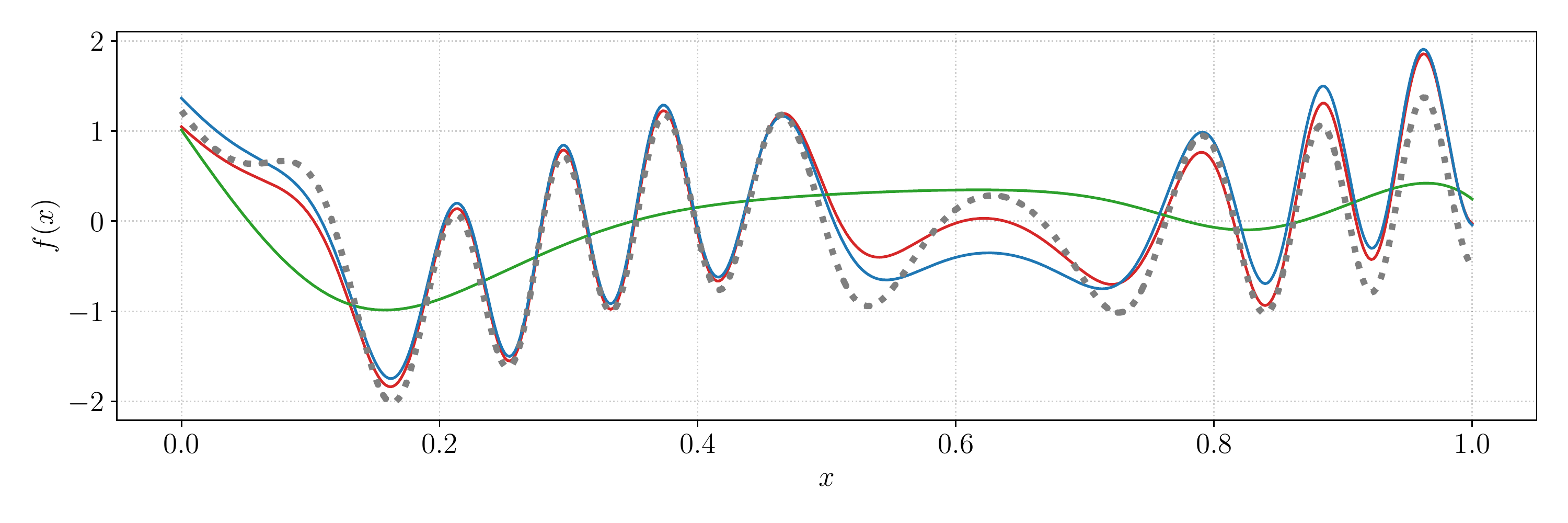}
    \caption{Predictions at two different test parameters using gappy POD, DeepONet, and multi-fidelity DeepONet with different configurations, varying the number of epochs and the POD energy threshold. The left column shows the results for $\mu = [4, 8]$, while in the right one we have $\mu = [3, 16]$. \emph{Top row}: POD energy threshold equal to 0.99, $10000$ epochs. \emph{Center row}: POD energy threshold equal to 0.999, $50000$ epochs. \emph{Bottom row}: POD energy threshold equal to 0.9999, $50000$ epochs.}
    \label{fig:appendix_1}
\end{figure}

\begin{figure}[H]
    \centering
    \includegraphics[width=.995\textwidth,trim={3cm 1cm 2cm 1cm},clip]{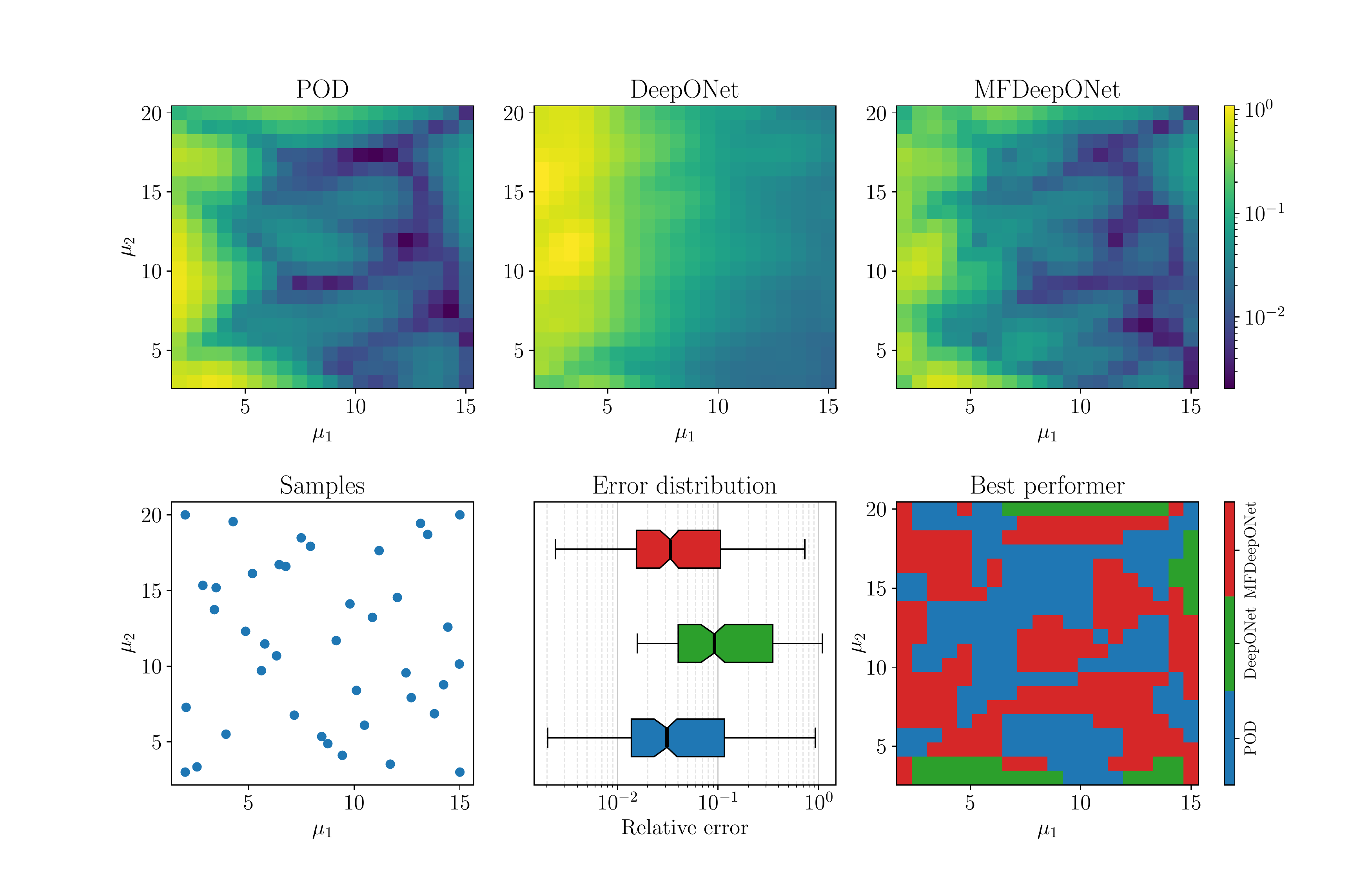}
    \caption{Comparison between POD ($0.9999$ energy threshold), DeepONet, and multi-fidelity
	DeepONet. From top to bottom we have the relative error in the
        parametric domain, the location of the high-fidelity samples in the
	parametric domain, the relative error distribution, and the best performers.}
    \label{fig:appendix_2}
\end{figure}

\begin{figure}[H]
    \centering
    \includegraphics[width=.99\textwidth,trim={.5cm .5cm .5cm 1cm},clip]{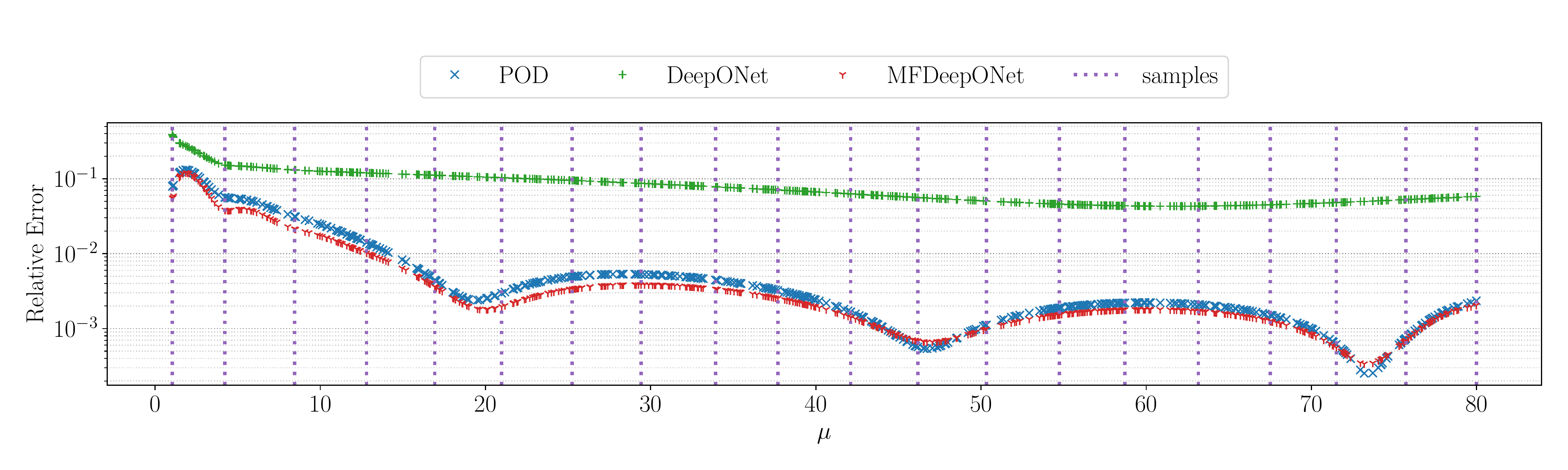}
    \caption{Comparison between POD (0.9999 energy threshold), DeepONet, and multi-fidelity
	DeepONet in terms of relative error in the parametric domain. The vertical dotted lines indicate the location of the high-fidelity samples.}
    \label{fig:appendix_3}
\end{figure}

This section presents additional plots for the POD energy threshold equal to $0.9999$ case, for the algebraic function (Figure~\ref{fig:appendix_1}) and for the parametric Navier-Stokes problem (Figures~\ref{fig:appendix_2} and~\ref{fig:appendix_3}).

\end{document}